%% file: main.tex
\documentclass[final,3p,times]{elsarticle}

\usepackage{amssymb}
\usepackage{amsmath}

\usepackage{bm}
\usepackage{soul}
\usepackage{float}

\usepackage{algorithm}
\usepackage{algpseudocode}

\usepackage{graphicx}
\usepackage{subcaption}

\usepackage{mathtools}

\usepackage{makecell}
\usepackage{hyperref}
\usepackage{xurl}
\usepackage{chngcntr}

\input{commands}

\journal{}

\begin{document}

\begin{frontmatter}

\title{Benchmark cases for gradient-based optimization in linear elastic solid mechanics using lattice Boltzmann methods}

\affiliation[mvm]{organization={Institute for Mechanical Process Engineering and Mechanics, Karlsruhe Institute of Technology},
             city={Karlsruhe},
             country={Germany}
}
\affiliation[ianm]{organization={Institute of Applied and Numerical Mathematics, Karlsruhe Institute of Technology},
             city={Karlsruhe},
             country={Germany}
}

\affiliation[lbrg]{
    organization={Lattice Boltzmann Research Group, Karlsruhe Institute of Technology},
     city={Karlsruhe},
     country={Germany}}

\affiliation[eth]{
    organization={Seminar for Applied Mathematics, ETH Z\"{u}rich},
    city={Zurich},
    country={Switzerland}} 

\author[mvm,lbrg]{Johannes L. Grafen}
\author[mvm,lbrg]{Florian Kaiser}
\author[ianm,lbrg,eth]{Stephan Simonis}
\author[mvm,lbrg]{Shota Ito}
\author[mvm,ianm,lbrg]{Mathias J. Krause}

\begin{abstract}
Gradient-based optimization (GBO) of problems constrained by partial differential
equations is central to parameter identification and structural design in solid
mechanics. Lattice Boltzmann method (LBM) schemes for linear elastic solids have only recently been derived. Combined with automatically generated discrete adjoint collision kernels, they make GBO possible with a purely LBM-based primal and adjoint solver.
We propose, to our knowledge, the first optimization benchmarks for solid LBM: three
inverse problems, each with a manufactured solution of the forward problem and a known reference control. These are the identification of two scalar amplitudes of a body force field on a periodic domain, the distributed control of the same field with up to $12\,800$ degrees of freedom, and the identification of Young's modulus on an elliptic plate with Dirichlet and Neumann boundaries. The primal scheme converges with second order on the periodic domain and with first order on the elliptic plate. The recovered controls follow these rates, with experimental orders of convergence of $1.92$ and $1.06$, respectively. Gradients obtained via forward-mode automatic differentiation (f-AD) and the discrete adjoint LBM (ALBM) require the same number of optimization steps and yield control errors agreeing to at least eleven significant digits. The ALBM gradient is further validated node-wise against a finite difference quotient, with deviations decreasing under mesh refinement. The benchmarks are implemented in the open-source LBM framework OpenLB and provide a verification basis for optimization in LBM-based solid mechanics and its extension to fluid-structure interaction.
\end{abstract}


\end{frontmatter}

\section{Introduction} \label{sec:intro}

Inverse problems and optimal control problems constrained by partial differential equations (PDEs) are a cornerstone of modern computational engineering. Solved by gradient-based optimization (GBO), they drive critical innovations in computational fluid and solid mechanics, including parameter identification \cite{pourasgharComputationallyEfficientApproach2024, bukshtynovVariationalStateDependentInverse2026}, shape reconstruction \cite{klemensCFDMRICoupledMeasurement2018, klemensSolvingFluidFlow2020, itoGeometryReconstructionMagnetic2026, cheylanShapeOptimizationUsing2019a} and topology optimization \cite{jenkinsImmersedBoundaryApproach2016, yoonTopologyOptimizationStationary2010, kuciLevelSetTopology2022, pingenTopologyOptimizationFlow2007a}, with applications ranging from aeroelasticity and acoustic-structure interaction to biomedical tissue scaffolding and the design of architected metamaterials \cite{pourasgharComputationallyEfficientApproach2024, cheylanShapeOptimizationUsing2019a, yoonTopologyOptimizationStationary2010, vergnaultAdjointbasedLatticeBoltzmann2014a, mauteConceptualDesignAeroelastic2004,munkTopologyOptimisationMicro2017,yoonBrittleDuctileFailure2017,kusanoAdjointbasedShapeOptimization2025}. What these problems have in common is that the optimizer requires the gradient of an objective functional with respect to a possibly high-dimensional control, evaluated subject to a PDE constraint that must be solved to (quasi-)steady state in every optimization step. 
In solid mechanics, where the constraint is traditionally discretized by the finite element method (FEM), the large matrix inversions required at high resolution remain a significant bottleneck for large-scale problems \cite{alexandersenReviewTopologyOptimisation2020}.

In the broader field of computational physics, the lattice Boltzmann method (LBM) has emerged as a highly efficient alternative to conventional numerical discretization schemes for the numerical solution of conservation equations \cite{simonisLatticeBoltzmannMethods2023}. 
LBM is exceptionally well-suited for large-scale parallelization on general-purpose graphics processing units (GPUs) \cite{kummerlanderImplicitPropagationDirectly2023, Krueger2016}. The OpenLB framework \cite{krauseOpenLBOpenSource2021, kummerlander_2025_17899765} has successfully scaled LBM simulations up to 512 GPUs on the HoreKa high-performance computing (HPC) cluster \cite{kummerlanderLargeScaleSimulationsTurbulent2025a}. To leverage these computational advantages for inverse problems, the continuous adjoint lattice Boltzmann method (ALBM) was originally proposed and implemented in OpenLB by \citet{krauseAdjointbasedFluidFlow2013}. Continuous ALBM formulations rely on a \textit{first-differentiate-then-discretize} approach, which requires the manual derivation of continuous adjoint equations for the specific LBM scheme \cite{klemensCFDMRICoupledMeasurement2018,klemensSolvingFluidFlow2020,krauseAdjointbasedFluidFlow2013, yajiTopologyOptimizationThermalfluid2016}. In contrast, discrete adjoint approaches follow the \textit{first-discretize-then-differentiate} paradigm. This guarantees consistency between the discretized objective functional and the computed gradient, and it allows the application of automatic differentiation (AD) software to differentiate arbitrary sets of equations \cite{gunzburgerPerspectivesFlowControl2002}. OpenLB \cite{krauseOpenLBOpenSource2021, kummerlander_2025_17899765} successfully applied the latter approach by combining AD and code generation with common subexpression elimination (CSE) to automatically generate discrete adjoint LBM collision kernels, resulting in speed-ups by a factor of three compared to manually derived continuous ALBM kernels \cite{itoGenerationEfficientAdjoint2026}.

Despite these advances, previous studies utilizing ALBM have predominantly focused on fluid dynamics, specifically Navier-Stokes flow control \cite{pingenTopologyOptimizationFlow2007a, krauseAdjointbasedFluidFlow2013, krauseParallelFluidFlow2013, laniewski-wollkAdjointLatticeBoltzmann2016}, thermal-fluid topology optimization \cite{yajiTopologyOptimizationThermalfluid2016, yajiTopologyOptimizationUsing2014, luoImprovedAdjointLattice2025b}, and advection-diffusion equations \cite{itoBenchmarkCaseInverse2025}. Even in complex multiphysics and fluid-structure interaction (FSI) optimization problems, researchers have historically restricted LBM to the modeling of the fluid phase while modeling the solid phase with traditional FEM solvers.

\citet{munkTopologyOptimisationMicro2017, munkEffectFluidstructureInteractions2018}, for example, performed multidisciplinary topology optimization by coupling LBM for the fluid phase with a fixed-mesh linear static FEM for the solid phase. Utilizing a bi-directional evolutionary structural optimization (BESO) algorithm, fluid pressures are mapped directly onto the structure as physical loads, the FEM is solved forward once, and sensitivities are extracted from the resulting elemental strain energies. 

\citet{cheylanDirectadjointLatticeBoltzmannSolver2026a} coupled an LBM flow solver with an FEM structural solver through the immersed boundary method (IBM) for turbulent FSI shape optimization. In their approach, the solid solution enters only in time-averaged form and is not taken into account in the sensitivity computation, so that the adjoint problem is restricted to the fluid phase. While effective, such couplings confine the computational advantages of LBM, and in this case also the gradient computation, to the fluid domain. Sensitivities with respect to the state of the solid itself are therefore not addressed in their work.

Recent foundational work in the derivation of LBM schemes for solid mechanics \cite{boolakeeNewLatticeBoltzmann2023, boolakeeDirichletNeumannBoundary2023}, the integration of these schemes into the OpenLB framework \cite{kaiserFluidStructureInteractionSimulations2025}, together with the automatic generation of efficient adjoint collision kernels \cite{itoGenerationEfficientAdjoint2026}, makes GBO with a purely LBM-based primal \emph{and} adjoint solver for solids possible for the first time. For optimization problems in LBM-based solid mechanics the literature currently offers neither a fully LBM-based optimization framework for solid mechanics nor rigorously validated benchmark cases for it.
The present work closes this gap. To the best of the authors' knowledge, it
constitutes the first application of the discrete ALBM to PDE-constrained
optimization problems in linear elastic solid mechanics in which the PDE
constraint is solved entirely by an LBM scheme. 

A benchmark case suited for this purpose has to meet several requirements. First, it must be reproducible, i.e., fully specified so that it can be reimplemented independently of a particular code. Second, it has to be
resolution-independent, the optimization problem and its solution have to be defined independently of the discretization, so that results obtained on different lattices converge to the same reference. Third, it has to allow for error quantification, both of the forward solution and of the recovered control, including their orders of convergence. Fourth, the computed gradient has to be verified against an independent reference. 
Within this work, we use manufactured solutions, where an analytical reference is being used as a reference state in the objective functional of the optimization problem.
The above requirements are addressed as follows. All geometries, material parameters, manufactured solutions, and optimizer settings are fully documented, and the cases are implemented in the OpenLB open-source framework \cite{krauseOpenLBOpenSource2021,
kummerlander_2025_17899765}. The overarching optimization problem is formulated in continuous space with a reference control $\al^{*}$ that does not depend on the lattice, and is solved over a sequence of resolutions. Manufactured solutions of the forward problem
give access to the primal discretization error. The known reference control $\al^{*}$ allows for the computation of the error of the recovered control, and the experimental order of convergence (EOC) of both is reported. Finally, the
gradients are verified against finite difference quotients (FDQ) and, where
feasible, by comparing two algorithmically independent differentiation methods, forward-mode AD (f-AD) and discrete ALBM. 

On this basis, we propose three benchmark cases derived from the simulation setups of \citet{boolakeeNewLatticeBoltzmann2023, boolakeeDirichletNeumannBoundary2023, kaiserFluidStructureInteractionSimulations2025}:
\begin{itemize}
  \item \textbf{Benchmark (1)} Parameter identification with periodic boundaries. The body force causing the deformation of a two-dimensional plate is determined inversely from two scalar control parameters, $\al \in \mathbb{R}^{2}$, for a component-wise scaling of an analytical force field. Gradients obtained from f-AD and from the discrete ALBM are compared.

  \item \textbf{Benchmark (2)} Distributed control problem with periodic boundaries. The full force field is inversely determined, such that the number of control variables scales with the number of nodes $\mathrm{dim}(\al) = 2N$. The gradients are computed with discrete ALBM and validated node-wise against forward FDQ.

  \item \textbf{Benchmark (3)} Parameter identification of Young's modulus on an elliptic plate with inner Dirichlet and outer Neumann boundaries using f-AD for gradient computation.
\end{itemize}
These cases are designed to serve as a standardized basis for GBO in LBM-based solid mechanics, ultimately paving the way for applications to non-linear elastic problems and FSI optimization.

This work is structured as follows. In \refSec{sec:problem}, the continuous PDE-constrained optimization problem is formulated. \refSec{sec:method} describes the numerical methods, including the LBM scheme for the primal problem, the discrete counterpart of the optimization problem, f-AD and the discrete ALBM as methods for gradient computation, and the employed GBO algorithm. The benchmarks and the obtained results are presented and discussed in \refSec{sec:experiments}. A summary and conclusion are given in \refSec{sec:conclusion}.

\section{Problem formulation} \label{sec:problem}
First, we define the continuous PDE-constrained optimization problem. Two types
of spatial domains are considered. For the bounded-domain benchmark,
$\Omega \subset \R^2$ is a bounded Lipschitz domain, with boundary $\partial \Omega = \partial \Omega_D \cup \partial \Omega_N$,
$\partial \Omega_D \cap \partial \Omega_N = \emptyset$, and outward unit normal
$\normal$. For the periodic benchmarks, $\Omega = \mathbb{T}^2 \coloneqq
\R^2/\mathbb{Z}^2$, so that $\partial \Omega = \emptyset$. In both cases, $\x \in \Omega$ denotes the position vecto. We define the steady-state displacement
$\disp \in \CV \coloneqq H^1(\Omega; \R^{2})$, where $\CV$ denotes the space of
admissible displacements, and the target displacement $\disp^* \in L^2(\Omega; \R^{2})$.

The solid material behavior is characterized by Young's modulus $E > 0$ and
Poisson's ratio $\nu$. Assuming a 2D plane stress state, the admissible range of
Poisson's ratio is $-1 < \nu < 1$, and the Lam\'e parameters relate to the material
parameters as follows
\begin{align}
    \lambda &= \frac{E \nu}{1 - \nu^2}, \label{eq:lambda} \\
    \mu &= \frac{E}{2(1+\nu)}, \label{eq:mu} \\
    K &= \mu + \lambda = \frac{E}{2(1-\nu)}, \label{eq:K}
\end{align}
where $K$ is referred to as the two-dimensional bulk modulus and
$\mu$ as the shear modulus \cite{boolakeeNewLatticeBoltzmann2023}. The associated
Cauchy stress tensor reads
\begin{equation}
    \stress(\disp) = K (\bm{\nabla} \cdot \disp) \mathbf{I}
    + \mu \Big[ \bm{\nabla} \disp + (\bm{\nabla} \disp)^{\mathsf{T}}
    - (\bm{\nabla} \cdot \disp) \mathbf{I} \Big],
    \label{eq:stress}
\end{equation}
and the applied body force is denoted by $\Force \in L^2(\Omega; \R^{2})$.

The admissible control space $\CS$ depends on the type of problem. Within this
work, we consider two types of problems
\begin{enumerate}
    \item Parameter identification problems, with
          $\CS \coloneqq \CSparam \subseteq \R^k$, $k \in \mathbb{N}$, that is, a
          finite-dimensional control space which is independent of the subsequent
          discretization.
    \item Distributed control problems, where the control is a square-integrable
          vector field that parameterizes the body force, either as the force
          itself or as a pointwise scaling of a given force field, such that
          $\CS \coloneqq \CSdist \subseteq L^2(\Omega;\R^2)$.
\end{enumerate}
The objective functional $\CJ: \CS \times L^2(\Omega; \R^{2}) \rightarrow \R$ measures the difference between the displacement $\disp$ and the target $\disp^*$ on the
observation domain $\Cdesign \subseteq \Omega$,
\begin{equation}
    \CJ(\al, \disp(\al) )\coloneqq \frac{1}{2}
    \displaystyle\int_{\Cdesign} \big\| \disp(\al, \x) - \disp^*(\x) \big\|_2^2 \, \mathrm{d}\x
    = \frac{1}{2} \big\| \disp(\al) - \disp^* \big\|_{L^2(\Cdesign)}^2 ,
    \label{eq:Jhat}
\end{equation}
and the optimization problem reads
\begin{equation}
    \min_{\al \in \CS}~\CJ(\al, \disp)
    \qquad \text{subject to} \qquad
    \hat{\bm{G}}(\al, \disp) = \bm{0},
    \label{eq:optproblem}
\end{equation}
where the constraint consists of the steady-state Navier--Cauchy equations of linear elasticity together with the boundary conditions,
\begin{equation}
    \hat{\bm{G}}(\al, \disp) \coloneqq
    \begin{cases}
        \mu\bm{\nabla}^2 \disp(\x) + K \bm{\nabla}(\bm{\nabla} \cdot \disp(\x)) + \Force(\x) = \bm{0}
            & \text{in} \quad \Omega, \\[4pt]
        \disp(\x) = \disp_D(\x)
            & \text{on} \quad \partial \Omega_D, \\[4pt]
        \stress(\disp(\x)) \, \normal = \trac(\x)
            & \text{on} \quad \partial \Omega_N,
    \end{cases}
    \label{eq:Ghat}
\end{equation}
with the prescribed Dirichlet displacement $\disp_D$ and the surface traction
$\trac$. The control enters \refEq{eq:Ghat} either through the body force,
$\Force = \Force(\al, \x)$, or through the material parameters $\mu(\al)$ and
$K(\al)$ via $E = E(\al)$; the specific dependencies are given with the respective
benchmarks in \refSec{sec:experiments}.

On the periodic domain, a steady state exists only if the body force satisfies the compatibility condition $\int_\Omega \Force \, \mathrm{d}\x = \bm{0}$, which restricts
the admissible controls in $\CS$ for the force-controlled benchmarks. The
displacement is then unique up to a rigid translation, consistent with the
zero-mean manufactured solution. We assume throughout that for every
$\al \in \CS$ the constraint admits a unique solution $\disp(\al) \in \CV$ that
depends differentiably on $\al$, so that the total derivative
$\mathrm{d}\CJ/\mathrm{d}\al$ used in \refSec{sec:method} is well defined.

Note finally that \refEq{eq:Ghat} is not solved directly. As described in
\refSec{sec:LBM_bulk}, it is recovered as the steady-state limit of a damped,
time-dependent lattice Boltzmann scheme.

\section{Solution procedure} \label{sec:method}
This section begins with outlining the LBM solution operator based on the work of
\citet{boolakeeNewLatticeBoltzmann2023}, which is used to solve the PDE constraint
(the primal or forward problem) in \refEq{eq:Ghat}, in \refSec{sec:LBM_bulk}. After
the description of the numerical method for the primal problem, the general
discrete inverse problem is formulated in \refSec{sec:discrete_problem}.
Subsequently, the two methods used for gradient computation in GBO are explained,
namely f-AD (\refSec{sec:AD}) and the discrete ALBM (\refSec{sec:ALBM}). Finally, the GBO algorithm is presented in \refSec{sec:GBO}.
\subsection{LBM scheme for solids} \label{sec:LBM_bulk}
In practice, the steady-state constraint \refEq{eq:Ghat} is solved using a
time-dependent LBM scheme. Following the approach of
\citet{boolakeeNewLatticeBoltzmann2023, boolakeeDirichletNeumannBoundary2023}, the
steady-state linear elasticity equation is extended by a time-dependent damping
term
\begin{equation}
    \kappa \partial_t \disp(\x,t) = \mu \bm{\nabla}^2 \disp(\x,t) + K \bm{\nabla}(\bm{\nabla} \cdot \disp(\x,t)) + \Force(\x) \quad \text{in} \quad \Omega \times [0,T] \subset \R^2 \times \mathbb{R}_0^+ ,
    \label{eq:NCE_t}
\end{equation}
where $\kappa$ denotes the damping constant and $T$ is the final time. In the
following, we set $\kappa=1$. Assuming that the body force $\Force(\x)$ is constant
in time, the solution of \refEq{eq:NCE_t} converges to the steady-state solution of
\refEq{eq:Ghat} for $T \rightarrow \infty$ \cite{boolakeeNewLatticeBoltzmann2023},
starting from an initial displacement $\disp(\x, t=0) = \disp_0(\x)$. The converged
field is then used as the state solution of the optimization problem.

\subsubsection*{Non-dimensionalization}

The application of the LBM requires a non-dimensionalization of the governing
equations. This is achieved by scaling physical quantities to their corresponding
lattice units, denoted by a tilde symbol \( (\tilde{\bullet}) \), using the
following relations
\begin{align}
    \x = L \tilde{\x}, \quad
    t=T \tilde{t}, \quad
    \disp = D \tilde{\disp}, \quad
    E = \frac{L^2 \kappa}{T} \tilde{E}, \quad
    \mu = \frac{L^2 \kappa}{T} \tilde{\mu}, \quad
    K = \frac{L^2 \kappa}{T} \tilde{K}, \quad
    \Force = \frac{L \kappa}{T} \tilde{\Force}, \quad
    \stress = \frac{L D \kappa}{T} \tilde{\stress}.
    \label{eq:unit_conversion}
\end{align}
Substituting these relations in \refEq{eq:NCE_t} yields the dimensionless form of
the Navier--Cauchy equation
\begin{equation}
    \partial_{\Tilde{t}}\Tilde{\disp} = \Tilde{\mu} \Tilde{\bm{\bm{\nabla}}}^2\Tilde{\disp} + \Tilde{K}\Tilde{\bm{\bm{\nabla}}}(\Tilde{\bm{\bm{\nabla}}} \cdot \Tilde{\disp}) + L D^{-1}\Tilde{\Force},
    \label{eq:NCE_t_dim}
\end{equation}
and the dimensionless Cauchy stress tensor is obtained as
\begin{equation}
    \Tilde{\stress} = \Tilde{K}(\Tilde{\bm{\bm{\nabla}}} \cdot \Tilde{\disp})\mathbf{I}
    + \Tilde{\mu} \Big[ \Tilde{\bm{\bm{\nabla}}}\Tilde{\disp} + (\Tilde{\bm{\bm{\nabla}}}\Tilde{\disp})^{\mathsf{T}} - (\Tilde{\bm{\bm{\nabla}}} \cdot \Tilde{\disp})\mathbf{I} \Big].
    \label{eq:stress_dim}
\end{equation}
The lattice moduli, and with them the relaxation times of
\refTab{tab:relaxation_times}, follow from the physical moduli as $\tilde{E} = E \Delta t / (\kappa \Delta x^2)$ (analogously for $\tilde{\mu}$ and
$\tilde{K}$). The pseudo-time step is therefore coupled to the grid spacing by the
diffusive scaling $\Delta t \propto \Delta x^2$, which keeps $\tilde{E}$ and the
relaxation times independent of the resolution. In all benchmarks of this work
$\Delta t = \Delta x^2$ is used, so that with $\kappa = 1$ the lattice moduli equal the
physical moduli, $\tilde{E} = E$. For the sake of readability, the tilde notation indicating dimensionless lattice units is dropped in the remainder of this work. If not stated otherwise, the reference quantities are chosen to be $1$.

\subsubsection*{Stencil and population vectors}
Using double-index notation, the multiple-relaxation-time (MRT) LBM scheme for
solving the time-dependent linear elasticity constraint \refEq{eq:NCE_t} is
formulated as in \cite{boolakeeNewLatticeBoltzmann2023, kaiserFluidStructureInteractionSimulations2025}.
The scheme uses a D2Q8 stencil with lattice speed $c = \Delta x / \Delta t$ and
discrete velocities $\bm{c}_{ij} = i c \, \bm{e}_x + j c \, \bm{e}_y$ for
$(i,j) \in \{-1, 0, 1\}^2 \setminus \{(0,0)\}$, notably omitting the rest population
$f_{00}$. The discrete populations are defined as
$f_{ij} = W_{ij} f(\x, t, \bm{c}_{ij})$ with the lattice weights $W_{ij}$. As in
\cite{boolakeeNewLatticeBoltzmann2023}, the weights enter only the definition of the
populations. The scheme itself is formulated entirely in terms of the $f_{ij}$, so
that their values are not required in the following.

Fixing the ordered set of lattice directions
\begin{equation}
    \mathcal{Q} \coloneqq \big( (1,0),\,(0,1),\,(-1,0),\,(0,-1),\,
                               (1,1),\,(-1,1),\,(-1,-1),\,(1,-1) \big),
    \label{eq:direction_order}
\end{equation}
indexed by $k \in \{1,\ldots,q\}$ with $(i_k,j_k) \in \mathcal{Q}$ and $q = 8$, the
local population vector of the node $\x_n$, $n \in \{0,\ldots,N-1\}$, reads
\begin{equation}
    \bm{f}_n(t) \coloneqq \big( f_{i_1 j_1}(\x_n,t), \ldots,
                                f_{i_q j_q}(\x_n,t) \big)^{\mathsf{T}} \in \R^{q},
    \label{eq:local_pop}
\end{equation}
where $N$ denotes the total number of lattice nodes. The post-streaming and
post-collision population vectors $\bm{f}_n^{\mathcal{S}}$ and $\bm{f}_n^{\mathcal{C}}$ are ordered analogously. For the global population vectors, the local ones are concatenated node by node,
\begin{equation}
    \fs(t) \coloneqq \big( \bm{f}_0^{\mathcal{S}}(t)^{\mathsf{T}}, \ldots,
                           \bm{f}_{N-1}^{\mathcal{S}}(t)^{\mathsf{T}}
                    \big)^{\mathsf{T}} \in \R^{qN},
    \label{eq:global_pop}
\end{equation}
and analogously for $\fc(t)$.

\subsubsection*{Moment space}
The ordering of $\mathcal{Q}$ matches the columns of the transformation matrix
$\mathbf{M} \in \R^{q \times q}$ given in \ref{app:transformation_matrices}, whose
rows yield the raw moments $m_{ab} = \sum_{(i,j)} i^a j^b f_{ij}$ in the order
\begin{equation}
    \bm{m} = \mathbf{M} \bm{f}_n^{\mathcal{S}}
    = \big( m_{10},\, m_{01},\, m_{11},\, m_s,\, m_d,\, m_{12},\, m_{21},\,
            m_{f} \big)^{\mathsf{T}} \in \R^{q},
    \label{eq:moments}
\end{equation}
with the combined higher-order moments
\begin{align}
    m_s &= m_{20} + m_{02}, \label{eq:moment_s}\\
    m_d &= m_{20} - m_{02}, \label{eq:moment_d}\\
    m_{f} &= m_{22} + \gamma\, m_s . \label{eq:moment_f}
\end{align}
The spherical/deviatoric split is due to \citet{boolakeeNewLatticeBoltzmann2023},
whereas the combination $m_f$ and the coefficient $\gamma$ were introduced by
\citet{boolakeeDirichletNeumannBoundary2023}, with
\begin{equation}
    \gamma = \frac{\theta \tau_f}{(1 + \theta)(\tau_s - \tau_f)},
    \qquad \theta = 1/3 .
    \label{eq:gamma}
\end{equation}
In the bulk $\gamma = 0$, so that $m_{f} = m_{22}$ and the last row of $\mathbf{M}$
reduces to the raw moment $m_{22}$. The coefficient $\gamma$ is only active in the
boundary collisions of the bounded-domain benchmark, where $\tau_f = 1/2$
(cf.\ \ref{app:transformation_matrices}).

The body force enters the scheme through the forcing term
\begin{equation} \label{eq:forcing}
    \bm{g}(\x) = \Delta t \, T^{-1} L D^{-1} \Force(\x),
\end{equation}
and is applied in two steps. The first half is added to the first-order moments,
which defines the bared first-order moments
\begin{align}
    \bar{m}_{10} &= m_{10} + g_x/2, \label{eq:moments_bar_1} \\
    \bar{m}_{01} &= m_{01} + g_y/2 . \label{eq:moments_bar_2}
\end{align}
The equilibrium moments are determined by these bared first-order moments alone,
\begin{equation}
    \bm{m}^{eq}
    = \mathbf{N}
      \begin{pmatrix} \bar{m}_{10} \\ \bar{m}_{01} \end{pmatrix},
    \qquad
    \mathbf{N}^{\mathsf{T}} =
    \begin{bmatrix}
        1 & 0 & 0 & 0 & 0 & \theta & 0      & 0 \\
        0 & 1 & 0 & 0 & 0 & 0      & \theta & 0
    \end{bmatrix}
    \in \R^{2 \times q}.
    \label{eq:moments_eq}
\end{equation}

\subsubsection*{Collision and streaming}
The post-collision moments follow from a moment-wise relaxation,
\begin{equation} \label{eq:moments_c}
    m^{\mathcal{C}}_{\beta} = m_\beta + \omega_{\beta} (m^{eq}_\beta - m_\beta),
    \quad \beta \in \mathcal{B} \coloneqq \{ 11,\, s,\, d,\, 12,\, 21,\, f \},
\end{equation}
where the relaxation time and the collision frequency are related by
\begin{equation} \label{eq:relaxation}
    \tau_\beta = 1/\omega_\beta - 1/2 ,
\end{equation}
with the moment-specific relaxation times listed in
\refTab{tab:relaxation_times}. The first-order moments are not relaxed, since
$m^{eq}_{10} = \bar{m}_{10}$ and $m^{eq}_{01} = \bar{m}_{01}$ by
\refEq{eq:moments_eq}. Instead, the second half of the forcing is applied to them,
\begin{align}
    m^{\mathcal{C}}_{10} &= \bar{m}_{10} + g_x/2, \label{eq:moments_c_10}\\
    m^{\mathcal{C}}_{01} &= \bar{m}_{01} + g_y/2. \label{eq:moments_c_01}
\end{align}
The bared moments are the arithmetic mean of their pre-and post-collision values,
\begin{equation} \label{eq:bared_moments}
    \bar{m}_{\beta} = (m_\beta + m_\beta^{\mathcal{C}}) / 2,
    \qquad \beta \in \mathcal{B},
\end{equation}
which is consistent with \refEq{eq:moments_bar_1}, since
$(m_{10} + m^{\mathcal{C}}_{10})/2 = m_{10} + g_x/2 = \bar{m}_{10}$. The
displacement and the Cauchy stress of a node are recovered from the bared moments as
\begin{align}
    \disp(\x_n) &= \begin{pmatrix} \bar{m}_{10} & \bar{m}_{01} \end{pmatrix}^{\mathsf{T}}(\x_n),
    \label{eq:disp_from_moments}\\[4pt]
    \stress(\x_n) &= -\begin{bmatrix}
        \bar{m}_{20} & \bar{m}_{11} \\
        \bar{m}_{11} & \bar{m}_{02}
    \end{bmatrix}(\x_n)
    = -\begin{bmatrix}
        (\bar{m}_s + \bar{m}_d)/2 & \bar{m}_{11} \\
        \bar{m}_{11} & (\bar{m}_s - \bar{m}_d)/2
    \end{bmatrix}(\x_n).
    \label{eq:stress_from_moments}
\end{align}
The collision and streaming steps of a single node can then be summarized as
\begin{align}
    \text{Collision:}& \quad
    f_{ij}^{\mathcal{C}}(\x, t) = \underbrace{
    \Big[ \mathbf{M}^{-1} \bm{m}^{\mathcal{C}}\big(\bm{f}_n^{\mathcal{S}}(\x, t), \Force(\x) \big) \Big]_{ij}
    }_{\eqqcolon\, \mathcal{C}_{ij}[
        \bm{f}_n^{\mathcal{S}}(\x, t), \Force(\x)
    ]}, \label{eq:local_collide} \\
    \text{Streaming:}& \quad
    f_{ij}^{\mathcal{S}}(\x + \bm{c}_{ij} \Delta t, t + \Delta t) = f_{ij}^{\mathcal{C}}(\x, t). \label{eq:local_stream}
\end{align}
Details about the implementation of this LBM scheme in OpenLB can be found in
\cite{kaiserFluidStructureInteractionSimulations2025}.

\subsubsection*{Operator notation}
In the following, we apply the commonly used operator notation for LBMs that fall
into the category of local lattice Boltzmann schemes (LLBS)
\cite{laniewski-wollkAdjointLatticeBoltzmann2016}, as in
\cite{itoGenerationEfficientAdjoint2026, laniewski-wollkAdjointLatticeBoltzmann2016, luoImprovedAdjointLattice2025b},
to denote the global collision and streaming updates of the global population
vectors \refEq{eq:global_pop},
\begin{align}
    \text{Collision:}& \quad
    \fc(t) = \mathcal{C}[ \fs(t), \Force],
    \label{eq:global_collide}\\
    \text{Streaming:}& \quad
    \fs(t + \Delta t) = \mathcal{S}[\fc(t)]
    \label{eq:global_stream}.
\end{align}
The global collision and streaming operators perform the mapping
$\mathcal{C}, \mathcal{S}: \R^{qN} \rightarrow \R^{qN}$. The operator
$\mathcal{C} = \mathrm{diag}(\mathcal{C}_n)$ has a block-diagonal matrix form, in
which the local collision operations \refEq{eq:local_collide} for the node with
index $n$ are contained in $\mathcal{C}_n$. The streaming operation is a linear
index-shifting operation
\cite{kummerlanderImplicitPropagationDirectly2023, laniewski-wollkAdjointLatticeBoltzmann2016},
such that the corresponding operator can be simplified to
$\mathcal{S}[\fc(t)] = \mathcal{S} \fc(t)$, where $\mathcal{S}$ is orthogonal
\cite{itoGenerationEfficientAdjoint2026}.

\subsection{Problem definition in discrete space} \label{sec:discrete_problem}
The continuous optimization problem defined in \refSec{sec:problem} is formulated in
the same discrete space in which the LBM scheme of the primal problem is solved. The
domain $\Omega$ is discretized by a uniform Cartesian lattice
$\Omega_{\Delta x} = \{ \x_0, \ldots, \x_{N-1} \}$ with spacing $\Delta x$ and node
index $n \in \{0, \ldots, N-1\}$ as in \refSec{sec:LBM_bulk}, where $N$ denotes the
total number of lattice nodes. For the square domains considered in this work,
$N = N_x^2$ with $N_x$ nodes per Cartesian direction, and $N_x$ serves as the
resolution parameter in the convergence studies of \refSec{sec:experiments}.
Accordingly, $\Ddesign \subseteq \Omega_{\Delta x}$ denotes the discretization of the
observation domain $\Cdesign$ on which the target displacement is available, and
$\mathcal{I}_{\Delta t} \coloneqq \{ 0, \Delta t, 2\Delta t, \ldots, T \}$ denotes the set of discrete time levels.

The discrete admissible control space $\CS_{\Delta x}$ is the discrete counterpart of
$\CS$ from \refSec{sec:problem}. For the parameter identification problems it
coincides with the continuous space, $\CS_{\Delta x} = \CSparam \subseteq \R^k$, since
a finite set of parameters does not depend on the discretization. For the distributed
control problem, the control field $\al \in \CSdist \subset L^2(\Omega; \R^2)$ is
represented by its nodal values, such that $\CS_{\Delta x} \subseteq \R^{2N}$ and the
dimension of the control scales with the resolution.

The state of the discrete system is the global post-streaming population vector
$\fs \in \R^{qN}$. The displacement field is recovered node-wise from the bared first-order moments \refEq{eq:disp_from_moments}, which depends on the local populations and, through the forcing term \refEq{eq:forcing}, additionally on the control itself. The discrete objective function is therefore a mapping
$\J: \CS_{\Delta x} \times \R^{qN} \rightarrow \R$ which evaluates the normalized
discrete $L^2$-error of the displacement field at the steady state $t=T$ of the
primal problem.

We first define a node-wise objective function
\begin{equation}
    \label{eq:node_J}
    J\big(\al, \bm{f}^{\mathcal{S}}_n\big) \coloneqq \frac{1}{2} \frac{\Big\| \disp\big(\al, \bm{f}^{\mathcal{S}}_n\big) - \disp^*(\x_n) \Big\|_2^2}{\bar{C}},
    \quad \text{with} \quad
    \bar{C} \coloneqq \sum_{\x_n \in \Ddesign} \Big\| \disp^*(\x_n) \Big\|_2^2 (\Delta x)^2 ,
\end{equation}
where $\bar{C}$ is the normalization constant. The sum of the node-wise evaluations of $J$ then yields the objective
functional $\J$ in discrete form, so that the discrete counterpart of the
optimization problem \refEq{eq:optproblem} reads
\begin{align}
    \min_{\al \in \CS_{\Delta x}} \J (\al, \fs) &= \sum_{\x_n \in \Ddesign} J\big(\al, \bm{f}^{\mathcal{S}}_n\big) (\Delta x)^2, \label{eq:discrete_J}\\
    \text{subject to} \quad
    \mathbf{G}(\al, \fs, \fc) &\coloneqq \begin{cases}
        \fc(t) - \mathcal{C}[\fs(t), \al], \\
        \fs(t + \Delta t) - \mathcal{S}\fc(t)
    \end{cases} = \bm{0}
    \quad \text{in} \quad \mathcal{I}_{\Delta t} \times \Omega_{\Delta x} , \label{eq:discrete_G}
\end{align}
closed by the initial populations $\fs(0)$ corresponding to the initial displacement
$\disp_0$.

Here, the continuous $L^2$-norm over $\Cdesign$ is approximated by a Riemann sum over
$\Ddesign$, so that $\bar{C}$ approximates $\| \disp^* \|^2_{L^2(\Cdesign)}$ and
$\J$ is the discrete counterpart of $\CJ$ from \refEq{eq:Jhat}, additionally
normalized in order to scale the objective function to a comparable range across
resolutions and test cases. Being a constant factor, this normalization does not
alter the location of the minimizer. The constraint \refEq{eq:discrete_G} corresponds
to the LBM scheme described in \refSec{sec:LBM_bulk}, written in operator notation
with the global population vectors $\fc$ and $\fs$
(\refEqs{eq:global_collide}{eq:global_stream}). Depending on the benchmark, the
control enters \refEq{eq:discrete_G} either through the forcing term
\refEq{eq:forcing}, $\Force = \Force(\al)$, or, for the identification of Young's
modulus $E$, through the lattice moduli and with them the relaxation frequencies
$\omega_\beta(\al)$ of the collision \refEq{eq:moments_c}.

The boundary conditions of the continuous problem \refEq{eq:Ghat} are not imposed as
separate constraints in \refEq{eq:discrete_G}. Instead, they are realized on the
mesoscopic level by modified collision and streaming operations at the boundary
nodes, and are hence already contained in the operators $\mathcal{C}$ and
$\mathcal{S}$ \cite{itoGenerationEfficientAdjoint2026}. Periodicity of the periodic benchmarks is recovered by the index shift
of the streaming operator, whereas the macroscopic Dirichlet and Neumann conditions of the bounded-domain benchmark are enforced by dedicated boundary collisions. For their derivation and implementation, the reader is referred to \cite{boolakeeDirichletNeumannBoundary2023, kaiserFluidStructureInteractionSimulations2025}.

\subsection{Gradient computation} \label{sec:gradients}
GBO methods require the total derivative $\dJda$ to perform the control update in
line-search based optimizers. Within the scope of this work, two methods are
utilized, f-AD and the discrete ALBM, which are explained in \refSec{sec:AD} and \refSec{sec:ALBM}, respectively.
\subsubsection{Automatic differentiation \label{sec:AD}}
Automatic differentiation (AD) evaluates derivatives of a computer program by systematically applying the chain rule to its elementary operations. In C++, AD is commonly realized by operator overloading, either by propagating or recording derivative information for every elementary operation at runtime
\cite{griewankAlgorithm755ADOLC1996}, or by expression templates, which assemble the derivative computation for entire expressions at compile time \cite{sagebaumHighPerformanceDerivativeComputations2019,
phippsEfficientExpressionTemplates2012}. Two modes of AD exist: the forward mode (f-AD), in which directional derivatives are propagated alongside the primal evaluation, and the reverse mode (r-AD), in which the elementary operations are recorded on a tape that is subsequently evaluated in reverse order to accumulate derivatives. OpenLB provides a f-AD type
based on operator overloading \cite{kummerlander_2025_17899765}, which is used in this work to compute gradients directly. r-AD is used only during the generation of the adjoint LBM collision kernels (\refSec{sec:ALBM}). For a
detailed explanation, the reader is referred to
\cite{itoGenerationEfficientAdjoint2026}. In the following, f-AD is described for the direct computation of the derivative of $\J$ with respect to the control variables $\al$.

In mathematical notation, f-AD corresponds to the following equation for a
differentiable function $\mathcal{F} : \R^n \rightarrow \R^m$ with
$\bm{x} \mapsto \bm{y} = \mathcal{F}(\bm{x})$ \cite{itoGenerationEfficientAdjoint2026}
\begin{equation} \label{eq:fAD}
    \dot{\bm{y}} = \mathbf{J}_{\mathcal{F}} (\bm{a}) \, \dot{\bm{x}},
\end{equation}
where the Jacobian of $\mathcal{F}$ is denoted by $\mathbf{J}_{\mathcal{F}}$ and is
evaluated at $\bm{a} \in \R^n$. Here, $\dot{\bm{x}} \in \R^n$ is the direction in
which the Jacobian is projected (also called the seeding vector) and
$\dot{\bm{y}} \in \R^m$ is the resulting directional derivative. Recovering the
full gradient therefore requires one primal evaluation per control variable. As
this renders f-AD computationally infeasible for distributed control problems
\cite{itoGenerationEfficientAdjoint2026}, f-AD is only used to compute $\dJda$ for
the parameter identification problems in \refTab{tab:overview_cases}.

\subsubsection{Adjoint LBM method \label{sec:ALBM}}
The adjoint LBM (ALBM) approach provides an efficient way to determine the gradient
$\dJda$ in distributed control problems, which is required by the GBO algorithm (see \refAlg{alg:GBO}). By the chain
rule, the total derivative of the objective function is
\begin{equation} \label{eq:dJda}
    \dJda = \fracpp{ \J }{\al} + \fracpp{\J}{\fs} \fracpp{\fs}{\al}.
\end{equation}
Explicitly computing the state sensitivity $\fracpp{\fs}{\al}$ is computationally
prohibitive for distributed control problems, where the control dimension satisfies
$\dim(\al) \gg 1$. To circumvent this, the Lagrangian function of the optimization
problem is introduced
\begin{equation} \label{eq:Lagrangian}
    \mathcal{L} = \J + \bm{\varphi}^{\mathsf{T}} \mathbf{G},
\end{equation}
where $\bm{\varphi} \in \R^{qN}$ is a vector of adjoint variables acting as Lagrange
multipliers, where $q$ denotes the number of populations per node and $N$ the total number of nodes. 

Because the governing equations are satisfied by the state solution
($\mathbf{G} = \bm{0}$), the total derivative of the objective function equals the total derivative of the Lagrangian for any choice of $\bm{\varphi}$,
\begin{equation} \label{eq:dLda_split}
    \dJda = \frac{d \mathcal{L}}{d \al}
    = \left[ \fracpp{\J}{\al} + \bm{\varphi}^{\mathsf{T}} \fracpp{\mathbf{G}}{\al} \right]
    + \left[ \fracpp{\J}{\fs} + \bm{\varphi}^{\mathsf{T}} \fracpp{\mathbf{G}}{\fs} \right]
      \fracpp{\fs}{\al},
\end{equation}
The term in the second bracket corresponds to $\fracpp{\mathcal{L}}{\fs}$ and setting it to zero yields the adjoint equation. 
This way, the state sensitivity $\fracpp{\fs}{\al}$ is eliminated without ever being formed, and the gradient reduces to
\begin{equation} \label{eq:dLda}
    \dJda = \frac{d\mathcal{L}}{d\al}
    = \fracpp{\J}{\al} + \bm{\varphi}^{\mathsf{T}} \fracpp{\mathbf{G}}{\al}.
\end{equation}
The ALBM follows a \textit{first-discretize-then-differentiate} approach
\cite{gunzburgerPerspectivesFlowControl2002}. In the following we refer to this method as discrete ALBM. Based on this approach, the adjoint-based optimization framework is proposed in \cite{itoGenerationEfficientAdjoint2026} and implemented in OpenLB \cite{krauseOpenLBOpenSource2021, kummerlander_2025_17899765}.

Since all benchmarks considered here are steady-state problems, the objective is
evaluated on the converged primal populations
$\fs^\infty := \lim_{t\to\infty} \fs(t)$ and $\fc^\infty := \lim_{t\to\infty} \fc(t)$.
Eliminating the time argument from \refEq{eq:discrete_G}, they satisfy the
fixed-point form of the constraint
\begin{equation} \label{eq:primal_fixed_point}
    \fc^\infty = \mathcal{C}\big[ \fs^\infty, \al \big],
    \qquad
    \fs^\infty = \mathcal{S} \fc^\infty
    \quad \Leftrightarrow \quad
    \fs^\infty = \mathcal{S}\, \mathcal{C}\big[ \fs^\infty, \al \big].
\end{equation}
Building the Lagrangian \refEq{eq:Lagrangian} on \refEq{eq:primal_fixed_point},
with the adjoint populations $\fsa$ and $\fca$ acting as multipliers of the
collision and the streaming part, respectively, and requiring
$\fracpp{\mathcal{L}}{\fs} = \fracpp{\mathcal{L}}{\fc} = \bm{0}$ yields the
adjoint fixed-point system in the same collide-and-stream structure as the primal
scheme,
\begin{alignat}{2}
    \text{Adjoint streaming:} \quad &\fsa^\infty = \mathcal{S}^{-1} \fca^\infty, \label{eq:adjoint_streaming} \\
    \text{Adjoint collision:} \quad &\fca^\infty = \Bigg[ \dCdf \Bigg]^{\mathsf{T}} \fsa^\infty - \Bigg[ \dJdf \Bigg]^{\mathsf{T}}, \label{eq:adjoint_collision}
\end{alignat}
where both Jacobians are evaluated at the converged primal state $\fs^\infty$ and
are therefore constant. Since $\mathcal{S}$ is orthogonal,
$\mathcal{S}^{-1} = \mathcal{S}^{\mathsf{T}}$. The linear system
\refEqs{eq:adjoint_streaming}{eq:adjoint_collision} is solved by the fixed-point
iteration
\begin{equation} \label{eq:adjoint_iteration}
    \fca^{(k+1)} = \Bigg[ \dCdf \Bigg]^{\mathsf{T}} \mathcal{S}^{-1} \fca^{(k)} - \Bigg[ \dJdf \Bigg]^{\mathsf{T}},
    \qquad k = 0, 1, 2, \ldots,
\end{equation}
started from
\begin{equation} \label{eq:adjoint_init}
    \fca^{(0)} = \bm{0},
\end{equation}
until the convergence criterion \refEq{eq:convergence_crit} is met, analogously to
the primal problem. This yields $\fca^\infty$ and
$\fsa^\infty$.

Evaluating \refEq{eq:dLda} for the fixed-point constraint
\refEq{eq:primal_fixed_point} gives the optimality condition
\begin{equation}
    \label{eq:optimality}
    \dJda = \Bigg[ \fracpp{\J(\al, \fs^\infty)}{\al} \Bigg]^{\mathsf{T}}
    - \Bigg[ \fracpp{\mathcal{C}[\fs^\infty, \al]}{\al} \Bigg]^{\mathsf{T}} \fsa^\infty .
\end{equation}
The control enters both terms only indirectly through an input variable, namely $\al \mapsto \Force(\al)$ for the force-controlled
benchmarks and $\al \mapsto E(\al) \mapsto \omega_\beta(E)$ for the identification of
Young's modulus. Writing $\J$ as a function of the populations
\refEq{eq:discrete_J} rather than of the displacement is what allows it to be
differentiated with respect to $\fs$ in \refEq{eq:adjoint_collision}. The explicit
dependence on $\al$ in \refEq{eq:disp_from_moments}, through the half-forcing
contribution of the bared moments, gives rise to the first term of
\refEq{eq:optimality}. The procedure is summarized in \refAlg{alg:ALBM}.
\begin{algorithm}[H]
\caption{Compute gradient $\dJda$ via Adjoint LBM (ALBM)}
\label{alg:ALBM}
\begin{algorithmic}[1]
\Require Converged primal state $\fs^\infty$, target displacement $\disp^*$, current control $\al^{(m)}$
\State Evaluate $\dJdf$, $\fracpp{\mathcal{C}[\fs^\infty,\al]}{\al}$ and $\fracpp{\J(\al,\fs^\infty)}{\al}$ \Comment{constant during the adjoint iteration}
\State Initialize $\fca^{(0)} = \bm{0}$, $k = 0$ \Comment{\refEq{eq:adjoint_init}}
\Loop { until adjoint convergence \refEq{eq:convergence_crit}}
    \State \textbf{Adjoint streaming:} \refEq{eq:adjoint_streaming}
    \State \textbf{Adjoint collision:} \refEq{eq:adjoint_collision}, $\Big[\dCdf\Big]^{\mathsf{T}}$ evaluated locally
    \State $k \gets k + 1$
\EndLoop
\State Set $\fca^\infty = \fca^{(k)}$, $\fsa^\infty = \mathcal{S}^{-1} \fca^\infty$ and evaluate \refEq{eq:optimality}
\State \Return $\left( \dJda \right)^{(m)}$
\end{algorithmic}
\end{algorithm}
The single-node Jacobian $\fracpp{\mathcal{C}}{\bm{f}^{\mathcal{S}}_n}$ has $q \times q$
entries, since the MRT collision \refEq{eq:moments_c} couples all populations of a
node through $\mathbf{M}$ and $\mathbf{M}^{-1}$. Storing it for the whole lattice would hence require $q^2 N$ floating-point numbers. That is $q$ times the
memory of the population field itself, and more than all remaining adjoint fields
combined. The ALBM framework \cite{itoGenerationEfficientAdjoint2026} therefore uses
AD to evaluate the Jacobians and partial derivatives of the adjoint system locally, and
combines code generation with common subexpression elimination (CSE), which removes the
runtime overhead of AD tapes and yields efficient adjoint kernels. The generated adjoint kernel evaluates the transposed Jacobian-vector product $\Big[ \fracpp{\mathcal{C}[\fs^\infty, \al]}{\bm{f}^{\mathcal{S}}_n} \Big]^T \fsai{n}^{(k)}$ locally in every adjoint
iteration (line 5 of \refAlg{alg:ALBM}), without ever assembling the global matrix. The remaining constant quantities require only $q$ entries per node and
are therefore evaluated once and stored. 

The generation of adjoint boundary kernels requires further extension to the existing code generation pipeline, as they include non-local operations \cite{boolakeeDirichletNeumannBoundary2023} and are therefore beyond the scope of this work. For further details on the derivation and implementation of the ALBM in OpenLB, the reader is referred to
\cite{itoGenerationEfficientAdjoint2026}.

\subsection{\label{sec:GBO} Gradient-based optimization algorithm}
The discrete optimization problem formulated in \refSec{sec:discrete_problem} is
solved iteratively using a gradient-based optimization (GBO) algorithm, conceptually
outlined in \refAlg{alg:GBO}. The algorithm starts with initial controls $\al^{(0)}$ and
iteratively updates the controls using a line search algorithm.

For a given optimization step $m$, the iteration consists of the following sequence.
First, the primal problem, which corresponds to the constraint
\refEq{eq:discrete_G}, is solved using the LBM scheme described in
\refSec{sec:LBM_bulk} until it reaches steady-state convergence. Convergence of both
the primal and the adjoint problem is determined by tracking the spatially averaged
squared norm of the displacement field $\disp$ (and similarly the adjoint displacement
field $\disp^{\varphi}$) over pseudo-time
\begin{equation} \label{eq:convergence}
    c(t) = \frac{1}{N} \sum_{n=0}^{N-1} \|\disp(\bm{f}^{\mathcal{S}}_n, t)\|_2^2.
\end{equation}
Steady state is considered reached when the coefficient of variation of $c(t)$ over
a trailing time window falls below a defined tolerance $\varepsilon_{conv}$
\begin{equation} \label{eq:convergence_crit}
    \frac{\mathrm{std}(c)}{\mathrm{mean}(c)} < \varepsilon_{conv},
\end{equation}
where $\mathrm{std}(c)$ and $\mathrm{mean}(c)$ denote the standard deviation and the mean of the
sequence, respectively.
It should be noted that $\varepsilon_{conv}$ can be chosen independently for the
primal and the adjoint problem.

Subsequently, the discrete objective functional \refEq{eq:discrete_J} is evaluated for the controls $\al^{(m)}$ in the current optimization step $m$. Next, the gradient of the objective functional
with respect to the control variables evaluated at $\al = \al^{(m)}$,
$\left( \dJda \right)^{(m)}$, is computed using either f-AD or the ALBM. The control parameters are then updated
using the quasi-Newton L-BFGS method \cite{liuLimitedMemoryBFGS1989}
\begin{equation} \label{eq:control_update}
    \al^{(m+1)} = \al^{(m)} + \delta \al^{(m)} \quad \text{with} \quad \delta \al^{(m)} = - \lambda^{(m)} \bm{p}^{(m)},
\end{equation}
where $\delta \al^{(m)}$ denotes the control update, $\lambda^{(m)}$ is the step size
calculated by the line search algorithm, and the descent direction
$\bm{p}^{(m)}$ is given by the L-BFGS approximation
\begin{equation} \label{eq:descent_direction}
    \mathbf{p}^{(m)} = \mathbf{H}^{(m)} \left( \dJda \right)^{(m)},
\end{equation}
with $\mathbf{H}^{(m)}$ representing the limited-memory approximation of the inverse
Hessian matrix. The step size $\lambda^{(m)}$ is determined via a line search that
must satisfy the strong Wolfe conditions
\cite{armijoMinimizationFunctionsHaving1966, wolfeConvergenceConditionsAscent1969, wolfeConvergenceConditionsAscent1971}.
An inner loop iteratively adapts the step size until these conditions are met
(compare \refAlg{alg:GBO}, lines 6 to 9). The optimization loop terminates if any of
the following criteria is satisfied: (1) the $L^2$-norm of the gradient falls below
a specified tolerance $\varepsilon$, (2) the relative change in the control
variables between successive steps is smaller than a threshold
$\varepsilon_{\alpha}$, or (3) the maximum allowed number of optimization steps
$m_{max}$ is reached.
\begin{algorithm}[H]
\caption{GBO algorithm with f-AD / ALBM and step-size control}
\label{alg:GBO}
\begin{algorithmic}[1]
\State Choose initial controls $\al^{(0)}$
\Loop { for \(m = 0,\,1,\,2,...,\,m_{\mathrm{max}}\)}
\State Solve primal problem $\mathbf{G}$ until steady-state convergence (\refEq{eq:convergence_crit}) and compute $\J(\al^{(m)}, \fs^\infty)$
\State Determine gradient $\left( \dJda \right)^{(m)}$ \Comment{via f-AD or ALBM (cf. \refAlg{alg:ALBM})}
\State Obtain descent direction $\bm{p}^{(m)}$ based on $\left( \dJda \right)^{(m)}$ \Comment{Using L-BFGS optimizer}
\Loop{ while $\lambda^{(m,k)}$ is invalid and $k<k_{max}$}
    \State{Choose new step size $\lambda^{(m,k)}$} \Comment{Step size control via strong Wolfe conditions}
    \State{Check if step is valid for $\delta \al^{(m,k)} = - \lambda^{(m,k)} \bm{p}^{(m)}$}
\EndLoop{}
\State Update controls $\al^{(m + 1)} = \al^{(m)} + \delta \al^{(m)}$ \Comment{Apply valid line search step}
\EndLoop { when $\Big|\Big| \left( \dJda \right)^{(m)} \Big|\Big|_2 < \varepsilon$
\textbf{or} for all $\alpha_i$ holds $\Big| 2\frac{\delta \alpha^{(m)}_i}{\alpha^{(m+1)}_i + \alpha^{(m)}_i} \Big| < \varepsilon_\alpha$} \Comment{Termination criteria}
\end{algorithmic}
\end{algorithm}

\section{Numerical experiments for solid LBM} \label{sec:experiments}
The three optimization benchmarks proposed in this work are summarized in
\refTab{tab:overview_cases}. All of them instantiate the general discrete problem
\refEqs{eq:discrete_J}{eq:discrete_G} and are built on primal structural mechanics
problems that admit an analytical (manufactured) solution, as first formulated in
\cite{boolakeeNewLatticeBoltzmann2023, boolakeeDirichletNeumannBoundary2023} and
implemented in OpenLB by \citet{kaiserFluidStructureInteractionSimulations2025}.
\begin{table}[H]
    \centering
    \begin{tabular}{l|c|c|c|c}
         & Method & Boundary conditions & Control  & Section \\
         \hline \hline
         (1) Parameter identification &  f-AD \& ALBM & periodic & \makecell{$\Force(\x_n, \al) = \al \circ \Force^*(\x_n)$,\\ with $\al \in \CS_{\Delta x} \subset \R^2$}& \ref{sec:parameter_periodic} \\
         \hline
         (2) Distributed control & ALBM & periodic & \makecell{$\Force(\x_n, \al(\x_n)) = \al(\x_n) = \al_i$,\\ with $\al \in \CS_{\Delta x} \subset \R^{2N}$} & \ref{sec:distributed_periodic} \\
         \hline
         (3) Parameter identification & f-AD & Neumann \& Dirichlet & \makecell{$E = \alpha^2$,\\ with $\alpha \in \CS_{\Delta x} \subset \R$} & \ref{sec:parameter_ellipse}\\
         \hline
    \end{tabular}
    \caption{Overview of the proposed optimization benchmarks.}
    \label{tab:overview_cases}
\end{table}
Benchmarks~(1) and~(2) share the periodic setup of
\cite{boolakeeNewLatticeBoltzmann2023} and differ in the dimension of the control
space, benchmark~(3) uses the bounded elliptic domain with Dirichlet and Neumann
boundaries of \cite{boolakeeDirichletNeumannBoundary2023}.

Each benchmark is investigated following the same three steps. First, the primal problem is verified against its manufactured solution in a grid convergence study, to determine the experimental order of convergence of the primal scheme. Second, the inverse problem is solved for a sequence of resolutions, and the error of the recovered control is measured against the known reference control $\al^*$. Third, the influence of the initial control $\al^{(0)}$ on the behavior of the GBO algorithm is assessed. For benchmark~(2) this procedure is complemented by a node-wise comparison of the ALBM gradient against a one-sided FDQ stencil.

\refSec{sec:numerical_validation} introduces the error measures used for the primal
verification. The benchmarks with periodic boundaries follow in
\refSec{sec:force_control}, subdivided into the numerical setup
(\refSec{sec:numerical_setup_periodic}), the verification of the primal problem
(\refSec{sec:primal_validation_periodic}), the parameter identification problem~(1)
(\refSec{sec:parameter_periodic}), the distributed control problem~(2)
(\refSec{sec:distributed_periodic}) and the gradient validation against the FDQ
stencil (\refSec{sec:gradient_validation}). The bounded-domain benchmark~(3) is
treated in \refSec{sec:param_bounded}, again subdivided into setup
(\refSec{sec:numerical_setup_ellipse}), primal verification
(\refSec{sec:primal_validation_ellipse}) and results
(\refSec{sec:parameter_ellipse}).

\subsection{Error measures and order of convergence} \label{sec:numerical_validation}
When analyzing the primal LBM scheme (\refSec{sec:LBM_bulk}) for linear elasticity,
the displacement and stress fields are of particular interest, since they constitute
the primary descriptors of solid structures subjected to external forces
\cite{gouldIntroductionLinearElasticity2018}. To assess the deviation from the
analytical solution, the numerical error of the primal problem is computed following
\cite{boolakeeNewLatticeBoltzmann2023}. 

At every lattice node, the numerical displacement and stress at the converged
steady state $t = T$ are reconstructed from the bared moments \refEqs{eq:disp_from_moments}{eq:stress_from_moments}
and are compared against the analytical solution $(\disp^*, \stress^*)$ introduced in the respective test cases. For a discrete function
$\chi: \Omega_{\Delta x} \rightarrow \R^{n}$ we employ the discrete norms
\begin{align}
    \| \chi \|_{L^2(\Omega_{\Delta x})} &= \left( (\Delta x)^2 \sum_{\x_n \in \Omega_{\Delta x}} \big| \chi(\x_n) \big|_2^{2} \right)^{\frac{1}{2}}, \label{eq:L2_norm}\\[6pt]
    \| \chi \|_{L^\infty(\Omega_{\Delta x})} &= \max_{\x_n \in \Omega_{\Delta x}} \big| \chi(\x_n) \big|_\infty, \label{eq:Linf_norm}
\end{align}
where $|\bullet|_2$ and $|\bullet|_\infty$ denote the Euclidean (Frobenius for
$\stress$) and the maximum norm of the nodal value, respectively. With
$p \in \{2, \infty\}$, the error of the primal problem is measured relative to the
analytical solution,
\begin{align}
    e^{\disp}_{L^p} &= \frac{\big\| \disp^* - \disp^{\mathrm{num}} \big\|_{L^p(\Omega_{\Delta x})}}
                           {\big\| \disp^* \big\|_{L^p(\Omega_{\Delta x})}},
    \label{eq:err_disp}\\[6pt]
    e^{\stress}_{L^p} &= \frac{\big\| \stress^* - \stress^{\mathrm{num}} \big\|_{L^p(\Omega_{\Delta x})}}
                             {\big\| \stress^* \big\|_{L^p(\Omega_{\Delta x})}},
    \label{eq:err_stress}
\end{align}
both evaluated at the converged steady state of the primal problem
(\refEq{eq:convergence_crit}), if not stated otherwise. 

For the inverse problems, the accuracy of the recovered control is assessed in the
same manner. Let $\al^{(M)}$ denote the control after the final optimization step
$M$, and $\al^*$ the known reference control of the respective benchmark. The
relative control error is then
\begin{equation} \label{eq:err_control}
    e^{\al} = \frac{\big\| \al^* - \al^{(M)} \big\|}{\big\| \al^* \big\|},
\end{equation}
where $\|\bullet\|$ denotes the Euclidean norm on $\CSparam \subseteq \R^k$ for the
parameter identification problems and the discrete $L^2$-norm
\refEq{eq:L2_norm} over $\Omega_{\Delta x}$ for the distributed control problem, in
which the control is a discrete function. 

Convergence rates for all of the above error measures, generically denoted by $e$,
are reported as an experimental order of convergence (EOC), obtained as the negative
slope of a least-squares fit of $\log e$ over $\log N_x$ across the resolutions
$N_x^{(r)}$ of the respective study,
\begin{equation} \label{eq:EOC}
    \mathrm{EOC} \coloneqq -\frac{\sum_{r} \big(\log N_x^{(r)} - \overline{\log N_x}\big)\big(\log e^{(r)} - \overline{\log e}\big)}
                          {\sum_{r} \big(\log N_x^{(r)} - \overline{\log N_x}\big)^2},
\end{equation}
where $r$ indexes the resolutions, $N_x$ denotes the number of lattice nodes per
Cartesian direction and $\overline{(\bullet)}$ the arithmetic mean. The sign is
chosen such that a scheme of order $\gamma$ yields $\mathrm{EOC} = \gamma$, since
$\Delta x \propto N_x^{-1}$.
\subsection{Optimal force control with periodic boundaries \label{sec:force_control}}
Benchmarks~(1) and~(2) share the periodic plate setup of
\cite{boolakeeNewLatticeBoltzmann2023}. In both benchmarks, the body force is the control and the
target is the manufactured displacement field, so that the exact force field $\Force^*$ is
known. They differ only in the dimension of the control space, in~(1) two scalar
parameters scale the components of $\Force^*$, independently of the resolution, while
in~(2) the nodal force field itself is the control, with $2N$ degrees of freedom.
This makes it possible to compare the f-AD and the ALBM gradient directly in~(1), and
to test the ALBM in the regime it is designed for in~(2). After the numerical setup
(\refSec{sec:numerical_setup_periodic}) and the verification of the primal problem
(\refSec{sec:primal_validation_periodic}), both inverse problems are solved on a sequence of resolutions and for different initial controls
(\refSec{sec:parameter_periodic} and \refSec{sec:distributed_periodic}). The ALBM
gradient of~(2) is finally validated node-wise against an FDQ reference (\refSec{sec:gradient_validation}). 
\subsubsection{Numerical setup \label{sec:numerical_setup_periodic}}
The first two benchmarks are posed on the unit square $\Omega = [0,1]^2 \subset \R^2$
with periodic boundaries aligned with the Cartesian unit vectors. The domain is
discretized by $N_x$ lattice nodes per Cartesian direction, that is $N = N_x^2$
nodes in total, and the observation domain covers the entire lattice,
$\Ddesign = \Omega_{\Delta x}$.

Following \cite{boolakeeNewLatticeBoltzmann2023}, a manufactured steady-state
solution $(\disp^*, \stress^*)$ is prescribed as
\begin{align}
    \disp^*(\x) &= \begin{pmatrix}
        \delta_1 \cos(2\pi x_1) \sin(2\pi x_2) \\
        \delta_2 \sin(2\pi x_1) \cos(2\pi x_2)
    \end{pmatrix}, \label{eq:disp_man_periodic}
\end{align}
with $\delta_1 = 9 \cdot 10^{-4}$, $\delta_2 = 7 \cdot 10^{-4}$, $\nu = 0.8$ and
$E = 0.11$ as in \cite{boolakeeNewLatticeBoltzmann2023}, where $\stress^*$ can be obtained from $\disp^*$ via \refEq{eq:stress}. 
The corresponding body force follows from the steady Navier--Cauchy equation as
\begin{equation} \label{eq:force_man_periodic}
     \Force^* = -\mu \bm{\nabla}^2 \disp^* - K\bm{\nabla} ( \bm{\nabla} \cdot \disp^*) = 
      4\pi^2 \begin{pmatrix}
        \big[K(\delta_1+\delta_2) + 2\mu\delta_1\big]\cos(2\pi x_1)\sin(2\pi x_2) \\[2pt]
        \big[K(\delta_1+\delta_2) + 2\mu\delta_2\big]\sin(2\pi x_1)\cos(2\pi x_2)
    \end{pmatrix},
\end{equation}
which satisfies the compatibility condition $\int_\Omega \Force^* \, \mathrm{d}\x = \bm{0}$
required in \refSec{sec:problem} and the residual rigid translation is fixed by the
zero-mean initial condition $\disp_0(\x, t=0) = \bm{0}$, which is imposed at the
beginning of every optimization step. The pair $(\disp^*, \Force^*)$ defines the
reference state of the inverse problems, indicated by the superscript
$(\cdot)^*$ (cf.\ \refEq{eq:Jhat}). The target displacement entering the objective
\refEq{eq:node_J} is the analytical field \refEq{eq:disp_man_periodic} sampled at
the lattice nodes.

\subsubsection{Verification of the primal problem \label{sec:primal_validation_periodic}}
\begin{figure}[!htbp]
    \centering
    \includegraphics{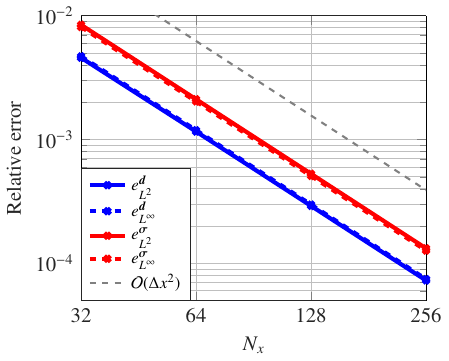}
    \caption{Grid convergence of the relative errors
    $e^{\disp}_{L^p}$ \refEq{eq:err_disp} and $e^{\stress}_{L^p}$ \refEq{eq:err_stress},
    $p \in \{2, \infty\}$, against the manufactured solution
    \refEqs{eq:disp_man_periodic}{eq:stress}, for
    $N_x \in \{32, 64, 128, 256\}$ lattice nodes per Cartesian direction, evaluated
    after a fixed number of $30 N_x^2$ time steps.}
    \label{fig:prim_periodic}
\end{figure}

To verify the implementation of the solid LBM scheme (\refSec{sec:LBM_bulk}), the
discretization error in displacement and stress is evaluated in a grid convergence
study for $N_x \in \{32, 64, 128, 256\}$. The results are shown in
\refFig{fig:prim_periodic}. To exclude the influence of the stopping criterion, every
run is performed with a fixed number of $30 N_x^2$ lattice steps ($T = 30$
with $\Delta t = \Delta x^2$) instead of the criterion \refEq{eq:convergence_crit}.
Second-order convergence is obtained in both $\disp$ and $\stress$ and in both norms,
with $\mathrm{EOC} = 1.99$ for $e^{\disp}_{L^2}$ and $e^{\disp}_{L^\infty}$ and
$\mathrm{EOC} = 2.00$ for $e^{\stress}_{L^2}$ and $e^{\stress}_{L^\infty}$
\refEq{eq:EOC}, in agreement with \cite{boolakeeNewLatticeBoltzmann2023, kaiserFluidStructureInteractionSimulations2025}.

\subsubsection{Parameter identification of the force field \label{sec:parameter_periodic}}
In the first benchmark the control $\al \in \CS_{\Delta x} \subset \R^2$ scales the
two components of the analytical force field \refEq{eq:force_man_periodic} evaluated
at the lattice nodes,
\begin{equation} \label{eq:control_param_periodic}
    \Force(\x_n, \al) = \al \circ \Force^*(\x_n)
    = \begin{pmatrix} \alpha_1 F^*_1(\x_n) \\ \alpha_2 F^*_2(\x_n) \end{pmatrix},
    \qquad n \in \{0,\ldots,N-1\}.
\end{equation}
By construction the reference control is
$\al^* = \begin{pmatrix} 1 & 1 \end{pmatrix}^{\mathsf{T}}$, and the exact minimizer
of the continuous problem is known. The control space is two-dimensional and independent of the resolution, so the same two-parameter problem is posed on every grid and only the discrete solution operator changes with $N_x$. The recovered control can therefore deviate from the reference control only through the error that the scheme itself commits on the manufactured solution. What the control error measures under grid refinement is thus the consistency error of the discretization, provided that the steady-state and optimizer tolerances are chosen small enough not to contribute, which is verified below. The parameters of the GBO algorithm and of the L-BFGS optimizer are listed in \refTab{tab:parameter}.
\renewcommand{\arraystretch}{1.75}
\begin{table}[H]
    \centering
    \begin{tabular}{|c|c|c|c|c|c|c|c|c|}
        \hline
        $T_{max}$ & $\varepsilon_{conv}$ & $\varepsilon$ & $\varepsilon_\alpha$ & $\lambda^0$ & $m_{max}$ & $k_{max}$ & $l$ & $\gamma^0$\\
        \hline
        $100$ & $10^{-7}$ & $10^{-10}$ &  $2.2 \cdot 10^{-16}$ & $1$ & $50$ & $100$ & $20$ & $10^{-4}$\\
        \hline
    \end{tabular}
    \caption{Parameters of the GBO algorithm and the L-BFGS optimizer for the
    parameter identification problem~(1), used identically for the f-AD and the ALBM
    gradient. $T_{max}$ is the maximum pseudo-time per primal and adjoint solve,
    $\varepsilon_{conv}$ the steady-state tolerance \refEq{eq:convergence_crit},
    $\varepsilon$ and $\varepsilon_\alpha$ the termination tolerances on the gradient
    norm and on the relative control change, $\lambda^0$ the initial line-search step size,
    $m_{max}$ and $k_{max}$ the maximum number of optimization and line-search
    iterations, $l$ the number of stored vector pairs of the limited-memory inverse
    Hessian approximation $\mathbf{H}^{(m)}$, and $\gamma^0$ the initial L-BFGS scaling.}
    \label{tab:parameter}
\end{table}
\renewcommand{\arraystretch}{1.0}
\begin{figure}[H]

    \begin{subfigure}{0.3\textwidth}
    \centering
    \includegraphics{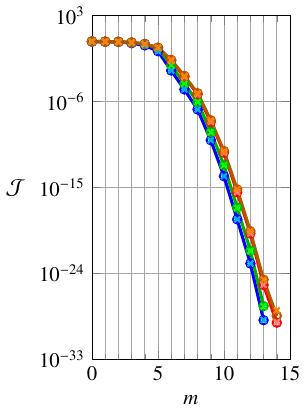}
    \label{fig:obj_param_AD}
    \end{subfigure}
    \hfill
    \begin{subfigure}{0.3\textwidth}
    \centering
    \includegraphics{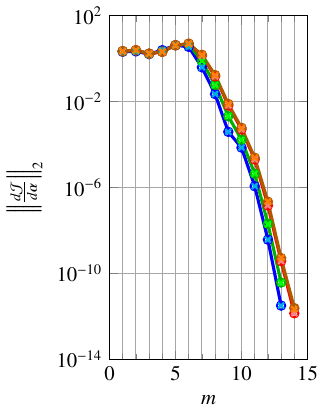}
    \label{fig:deriv_param_AD}
    \end{subfigure}
    \hfill
    \begin{subfigure}{0.3\textwidth}
    \centering
    \includegraphics{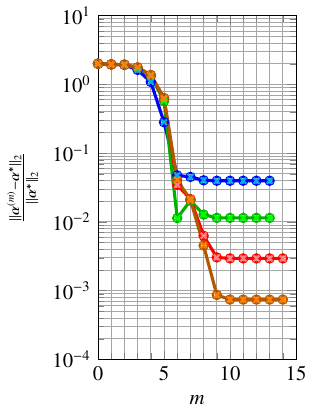}
    \label{fig:controls_param_AD}
    \end{subfigure}
    \vspace{5pt}
    \centering \includegraphics{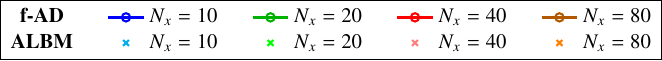}

    \caption{Objective functional $\J$, Euclidean norm of the gradient
    $\big\|\dJda\big\|_2$ and relative error of the controls
    $\relErrorParameterControl$ for the parameter identification problem~(1) of
    \refTab{tab:overview_cases}, with the gradient computed by f-AD and by the ALBM, for resolutions
    $N_x \in \{10,20,40,80\}$ and initial control
    $\al^{(0)} = \begin{pmatrix} -1 & -1 \end{pmatrix}^{\mathsf{T}}$.}
    \label{fig:param_periodic}
\end{figure}
\refFig{fig:param_periodic} shows the evolution of the objective functional $\J$,
of the gradient norm $\|\dJda\|_2$ and of the relative control error \refEq{eq:err_control} over the optimization steps $m$, for $N_x \in \{10,20,40,80\}$ and for both gradient computation methods, starting from
$\al^{(0)} = \begin{pmatrix} -1 & -1 \end{pmatrix}^{\mathsf{T}}$.

The f-AD and the ALBM results agree in all three quantities and at every resolution. Both methods require the same number of optimization steps, and the final control errors agree to at least eleven significant digits.
Since the two methods share only the
primal solution, their agreement is a non-trivial check of the generated adjoint collision kernel.

The optimization terminates within $15$ steps at all resolutions, and the
objective decreases to the level of double-precision round-off. The reason is that
the target is exactly representable by the two available control parameters. 
The scheme of \refSec{sec:LBM_bulk} is linear in the populations and in the body
force.

The relaxation frequencies depend only on the material parameters. Since, in addition, all coefficients are the same at every node of the uniform periodic lattice, the scheme does not generate spatial frequencies that are absent from the forcing. The body force \refEq{eq:force_man_periodic} contains a single wavenumber, so the converged numerical solution $\disp^{(k)}_{\Delta x}$ obtained for $\Force^{opt} = \Force^*$, retains this wavenumber and differs from the manufactured solution \refEq{eq:disp_man_periodic} only through the amplitudes of its two components, which carry the discretization error. The residual $\disp_{\Delta x}(\al) - \disp^*$ is therefore driven to round-off, and $\J$, being quadratic in that residual, stagnates between $10^{-30}$ and $10^{-27}$. 
The relative control error, in contrast, stagnates at a resolution-dependent floor, because the control that achieves the exact fit is not $\al^*$. At $\al^*$ the scheme reproduces the manufactured solution only up to the amplitude error of the discretization. The control that fits the target exactly therefore rescales the two force components by the inverse of that error, so that $\al^{opt} - \al^*$ is the relative amplitude error committed by the scheme, which is of $\mathcal{O}(\Delta x^2)$. The floor is thus a property of the discretization and not of the optimizer, which is confirmed quantitatively. The EOC of the relative control error in the final optimization step, measured over $N_x \in \{10,20,40,80\}$ according to \refEq{eq:EOC}, is $1.919$ for both f-AD and the ALBM. Since the manufactured solution consists of a single wavenumber, this amplitude error is also what the primal error of \refFig{fig:prim_periodic} measures, which explains why both quantities converge with second order. The relative control error therefore measures the consistency error of the discretization, and not a property of the gradient computation or of the optimizer. Convergence of the primal and adjoint solvers is monitored with the criterion \refEq{eq:convergence_crit} using $\varepsilon_{conv} = 10^{-7}$, where a simulation run is terminated at $T_{max}$ (see \refTab{tab:parameter}) if the criterion is not met before.

\begin{figure}[H]

    \begin{subfigure}{0.3\textwidth}
    \centering
    \includegraphics{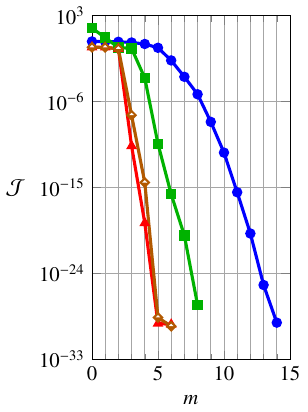}
    \label{fig:obj_param_AD_init}
    \end{subfigure}
    \hfill
    \begin{subfigure}{0.3\textwidth}
    \centering
    \includegraphics{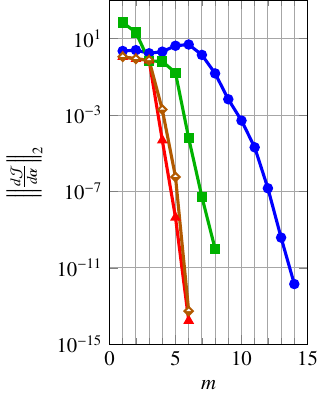}
    \label{fig:deriv_param_AD_init}
    \end{subfigure}
    \hfill
    \begin{subfigure}{0.3\textwidth}
    \centering
    \includegraphics{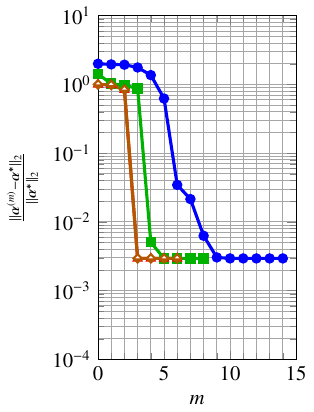}
    \label{fig:controls_param_AD_init}
    \end{subfigure}
    \vspace{5pt}
    \centering
    \includegraphics{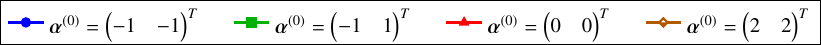}

    \caption{Objective functional $\J$, Euclidean norm of the gradient
    $\big\|\dJda\big\|_2$ and relative error of the controls
    $\relErrorParameterControl$ for the parameter identification problem~(1) of
    \refTab{tab:overview_cases}, computed with f-AD for the four initial controls
    $\al^{(0)} \in \big\{ (-1, -1)^{\mathsf{T}}, (-1, 1)^{\mathsf{T}},
    (0,0)^{\mathsf{T}}, (2,2)^{\mathsf{T}} \big\}$ at resolution $N_x = 40$.}
    \label{fig:param_init_controls_periodic}
\end{figure}

The influence of different initial controls is examined in
\refFig{fig:param_init_controls_periodic}, where four initial controls
\begin{equation*}
    \al^{(0)} \in \{ (-1,-1)^{\mathsf{T}}, (-1,1)^{\mathsf{T}}, (0,0)^{\mathsf{T}},
(2,2)^{\mathsf{T}} \}
\end{equation*}
are compared at $N_x = 40$ with the f-AD gradient. All four runs reach the same objective and control-error floor.

\subsubsection{Distributed force control} \label{sec:distributed_periodic}
\begin{figure}[!htbp]
    \centering
    \begin{subfigure}{0.3\textwidth}
        \centering
        \includegraphics{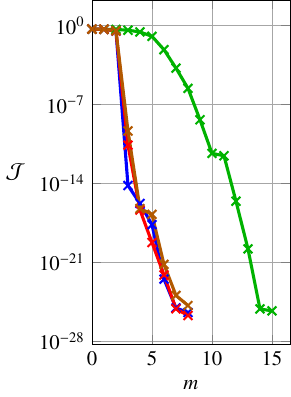}
        \label{fig:obj_dist_ALBM}
    \end{subfigure}
        \hspace{0.02\textwidth}
    \begin{subfigure}{0.3\textwidth}
        \centering
        \includegraphics{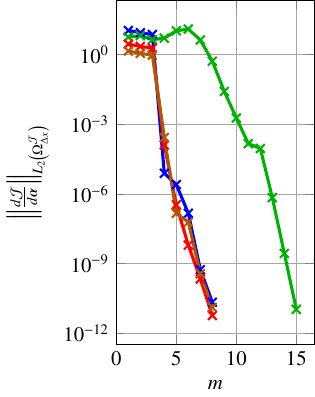}
        \label{fig:deriv_dist_ALBM}
    \end{subfigure}
        \hspace{0.02\textwidth}
    \begin{subfigure}{0.3\textwidth}
        \centering
        \includegraphics{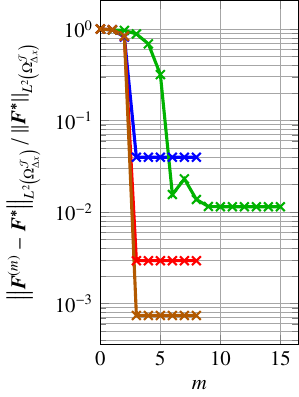}
        \label{fig:controls_dist_ALBM}
    \end{subfigure}
        \hfill
    \begin{subfigure}{0.3\textwidth}
        \centering
        \includegraphics{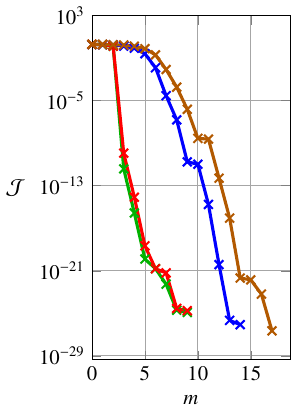}
        \label{fig:obj_dist_ALBM_informed}
    \end{subfigure}
        \hspace{0.02\textwidth}
    \begin{subfigure}{0.3\textwidth}
        \centering
        \includegraphics{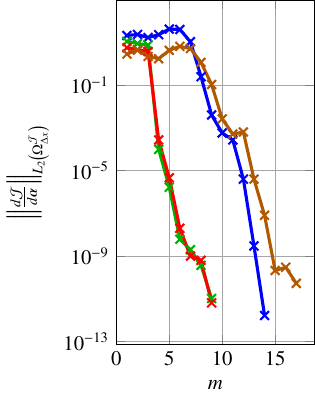}
        \label{fig:deriv_dist_ALBM_informed}
    \end{subfigure}
        \hspace{0.02\textwidth}
    \begin{subfigure}{0.3\textwidth}
        \centering
        \includegraphics{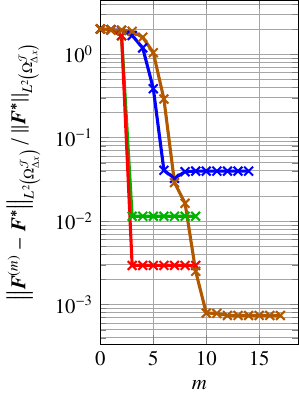}
        \label{fig:controls_dist_ALBM_informed}
    \end{subfigure}
    \vspace{5pt}
    \centering \includegraphics{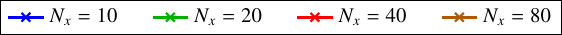}

    \caption{Objective functional $\J$, $L^2$-norm of the gradient $\normdispderivative$ and relative $L^2$-error of the force field $\relErrorControl{m}$ for the distributed control problem of
    \refTab{tab:overview_cases}, computed with the ALBM for resolutions
    $N_x \in \{10,20,40,80\}$, with the zero-mean projection
    \refEq{eq:zero_mean_projection}. Top row: uninformed initial control~(i),
    $\al^{(0)}(\x_n) = \begin{pmatrix} -1 & -1 \end{pmatrix}^{\mathsf{T}}$. Bottom row: informed initial control~(ii),
    $\al^{(0)} = -\Force^*$. The optimizer parameters are given in
    \refTab{tab:distributed}.}
    \label{fig:dist_ALBM_1}
\end{figure}

In the second benchmark the force field itself is the control, represented by its
nodal values, so that $\al \in \CS_{\Delta x} \subset \R^{2N}$ with
$\al = \al(\x_n)$, $n \in \{0, \ldots, N-1\}$. The number of degrees of freedom
$\mathrm{dim}(\al) = 2N$ therefore scales with the resolution, which renders f-AD
computationally infeasible (cf.\ \refSec{sec:AD}). All gradients in this section are computed with the ALBM. The control is the nodal force field,
$\Force(\x_n, \al) = \al(\x_n)$, with the reference control
$\al^*(\x_n) = \Force^*(\x_n)$. Two initial controls are considered,
\begin{enumerate}
    \item[(i)] the uninformed initial control
               $\al^{(0)}(\x_n) = \begin{pmatrix} -1 & -1 \end{pmatrix}^{\mathsf{T}}$, which
               carries no information about the analytical solution,
    \item[(ii)] the informed initial control $\al^{(0)}(\x_n) = -\Force^*(\x_n)$.
\end{enumerate}
On the periodic domain a static equilibrium exists only for body forces with
vanishing mean (cf.\ \refSec{sec:problem}). A control with non-zero mean, such as
the constant initial control of (i), violates this compatibility condition. This causes the
displacement to drift, $\J$ then depends on $T_{max}$ rather than on a steady state,
and the steady adjoint is no longer consistent with the primal solve. The control and the
gradient are therefore projected component-wise onto the zero-mean subspace,
\begin{equation} \label{eq:zero_mean_projection}
    \mathcal{P}\al(\x_n) = \al(\x_n) - \frac{1}{N} \sum_{k=0}^{N-1} \al(\x_k).
\end{equation}
With the projected initial control and projected gradients, the L-BFGS iterates remain in this subspace. For (i) the projection maps the constant initial control to
$\al^{(0)} = \bm{0}$, whereas the initial control of (ii) already has zero mean and is left
unchanged. The error of the recovered control is reported for the force field,
$\relErrorControl{m}$, and the gradient in the corresponding $L^2$-norm over
$\Ddesign$. The parameters of the GBO algorithm and of the L-BFGS optimizer are listed
in \refTab{tab:distributed}.
\renewcommand{\arraystretch}{1.75}
\begin{table}[H]
    \centering
    \begin{tabular}{|c|c|c|c|c|c|c|c|c|}
        \hline
        $T_{max}$ & $\varepsilon_{conv}$ & $\varepsilon$ & $\varepsilon_\alpha$ & $\lambda^0$ & $m_{max}$ & $k_{max}$ & $l$ & $\gamma^0$\\
        \hline
        $10$ & $10^{-7}$ & $10^{-10}$ &  $2.2 \cdot 10^{-16}$ & $1$ & $100$ & $100$ & $20$ & $10^{-4}$\\
        \hline
    \end{tabular}
    \caption{Parameters of the GBO algorithm and the L-BFGS optimizer for the
    distributed control problem~(2), with the symbols as defined in
    \refTab{tab:parameter}.}
    \label{tab:distributed}
\end{table}
\renewcommand{\arraystretch}{1.0}
\refFig{fig:dist_ALBM_1} shows the evolution of the objective functional $\J$, of the
gradient norm and of the relative error of the recovered force field over the
optimization steps $m$, for $N_x \in \{10,20,40,80\}$ and for both initial controls. In both cases the objective decreases from its initial value to $\J \approx \mathcal{O}(10^{-25})$, the level of double-precision round-off.
The relative error of the recovered force field, in contrast, again stagnates at a resolution-dependent floor.  Since the control has $2N$ degrees of freedom, the sampled target is attainable exactly, and the recovered field is the discrete force that reproduces it. Inserting the analytical solution into the discrete equations leaves a residual, the consistency error of the scheme, so that $\Force^{opt}$ differs from $\Force^*$ by exactly this residual and $|\Force^{opt} - \Force^*|_{L^2(\Ddesign)} = \mathcal{O}(\Delta x^2)$. The measured EOC of the relative force error in the final optimization step is $1.919$ for both
initial controls. The final errors are $3.97 \cdot 10^{-2}$, $1.14 \cdot 10^{-2}$, $2.94 \cdot 10^{-3}$ and $7.41 \cdot 10^{-4}$ for $N_x \in \{10, 20, 40, 80\}$ and agree to at
least ten digits between the two initial controls. 

The two initial controls differ in their starting point. After the projection, (i) starts
from $\al^{(0)} = \bm{0}$, that is from $\J = 1/2$ and a relative force error of $1$, and
requires at most $15$ optimization steps. (ii) starts from $\Force^{(0)} = -\Force^*$ at a relative error of $2$ and requires at most $17$ steps. Both reach the same floor. 

\subsubsection{Gradient validation \label{sec:gradient_validation}}
In this section the ALBM gradient of the
distributed control problem is validated node-wise against
a one-sided FDQ stencil, as in
\cite{itoGenerationEfficientAdjoint2026, tekitekAdjointLatticeBoltzmann2006a},
\begin{equation} \label{eq:FDQ_grad}
    \frac{d \J}{d\alpha_{i,j}} \Bigg|^{FDQ}_{\x = \x_n}
    = \frac{ \J(\al + h \bm{e}_{i,j}) - \J(\al)}{h} + \mathcal{O}(h) \quad \text{with} \quad h=10^{-5},
\end{equation}
where $\bm{e}_{i,j}$ denotes the unit vector of the control component
$j \in \{x,y\}$ at the node $\x_n$ and $h$ is the perturbation. The comparison of the gradients is performed for a subset of interior nodes along the diagonal $y = x$, parameterized by
$\x_n = \begin{pmatrix} s & s \end{pmatrix}^{\mathsf{T}}$ with
$s \in \{\Delta x, 2\Delta x, \ldots, 1 - \Delta x\}$, excluding the periodic
boundary at $s \in \{0,1\}$. For $N_x > 10$ the nodes of the $N_x = 20$ lattice are used, so that all resolutions are evaluated at the same physical positions. The gradient comparison is reported for $h=10^{-5}$, the largest step size for which all perturbed forward solves converge so that the deviations given below are upper bounds for the error of the ALBM gradient.

\begin{figure}[!htbp]
  \centering
    \begin{subfigure}[t]{0.33\textwidth}
    \centering
    \raisebox{-139.59348pt}{\includegraphics{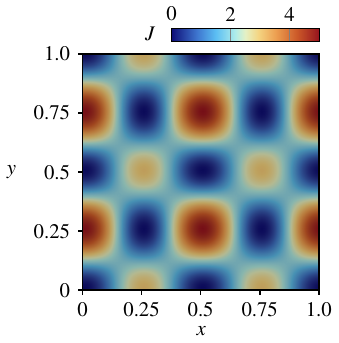}}
    \caption{Node-wise objective $J(\x)$}
    \label{fig:objective}
    \end{subfigure}
      \hfill
    \begin{subfigure}[t]{0.33\textwidth}
    \centering
    \raisebox{-139.52296pt}{\includegraphics{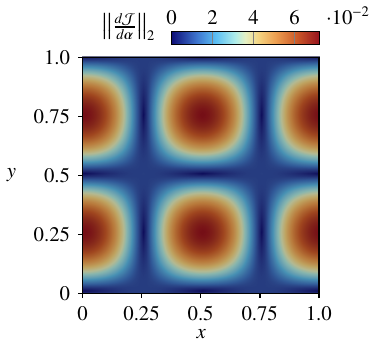}}
    \caption{Magnitude of $\dJda$}
    \label{fig:sensitivity}
    \end{subfigure}
    \hfill
    \begin{subfigure}[t]{0.33\textwidth}
        \centering
        \raisebox{-138.97107pt}{\includegraphics{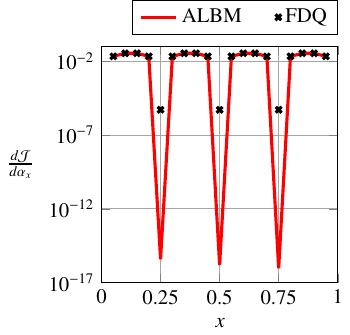}}

        \caption{Comparison of $\frac{d\J}{d \alpha_x} \Big |_{\bm{x} = \bm{x_i}}$ for ALBM and FDQ}
        \label{fig:grad_validation_res_80}
    \end{subfigure}

  \caption{Initial node-wise objective functional \refFig{fig:objective} and magnitude of sensitivities \refFig{fig:sensitivity} computed with the ALBM for resolution $N_x=80$ and initial force field $\Force^{(0)} = -\Force^*$. The respective sensitivities in x-direction are evaluated for a subset of points and compared to the FDQ gradient in \refFig{fig:grad_validation_res_80}.
  }
  \label{fig:objective_sensitivity}
\end{figure}

The gradients are validated for the informed initial control~(ii) of
\refSec{sec:distributed_periodic} in the first optimization step,
$\al^{(0)} = \Force^{(0)} = -\Force^*$. Since this control already has zero mean, the gradients are compared without the projection
\refEq{eq:zero_mean_projection}. The primal and adjoint problems are solved with the parameters given in \refTab{tab:distributed}. The node-wise objective $J(\x)$ and the magnitude of the corresponding sensitivity field obtained from \refEq{eq:optimality} are shown in
\refFig{fig:objective} and \refFig{fig:sensitivity} for
$N_x = 80$.
\refFig{fig:grad_validation_res_80} compares the two gradients along the diagonal.
They agree over the entire range except at three isolated locations, where the ALBM
gradient drops by up to fifteen orders of magnitude while the FDQ value stagnates.
These locations are not incidental. The control enters the scheme only through the
forcing \refEq{eq:forcing}, so $\partial \bm{g} / \partial \al$ is a constant
conversion factor and the gradient \refEq{eq:optimality} is proportional to the adjoint
displacement. For $\Force^{(0)} = -\Force^*$ the residual is
$\disp^{(0)} - \disp^* \approx -2\disp^*$, which contains the single wavenumber of
\refEq{eq:disp_man_periodic}. The linear adjoint scheme with constant coefficients
preserves this spatial structure, so that the $x$- and $y$-components of the gradient
are proportional to the respective components of $\disp^*$. Along the diagonal both components of $\dJda$ are proportional to $\tfrac{1}{2}\sin(4\pi s)$, whose
interior zeros lie at $s \in \{0.25, 0.5, 0.75\}$ and coincide with the observed
positions, where the exact gradient vanishes. At such a point the one-sided stencil \refEq{eq:FDQ_grad} evaluates to its truncation term $\mathcal{O}(h)$ rather than zero, and its stencil width
limits the attainable precision.

The absolute and relative deviations from the FDQ reference are presented in
\refFig{fig:error_grad} for the resolutions $N_x \in \{ 10,20,40,80\}$. The gradients are compared at the interior diagonal nodes, excluding the periodic boundary. The deviation of the ALBM from the FDQ gradients along the diagonal decreases for an increasing spatial resolution, which confirms the validity of the ALBM
gradient.
\newcommand{\dJdax}{\frac{d \J}{d \alpha_x}}
\newcommand{\dJday}{\frac{d \J}{d \alpha_y}}
\begin{figure}[!htbp]
    \begin{subfigure}{0.49\textwidth}
        \centering
        \includegraphics{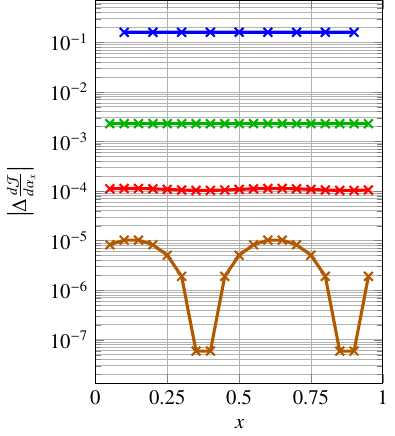}
        \caption{Absolute errors of $\dJdax$}
        \label{fig:abs_error_grad_x}
    \end{subfigure}
    \hfill
    \begin{subfigure}{0.49\textwidth}
        \centering
        \includegraphics{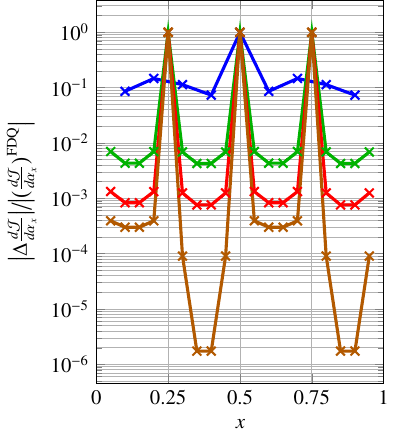}
        \caption{Relative errors of $\dJdax$}
        \label{fig:rel_error_grad_x}
    \end{subfigure}
    \par\vspace{0.8cm}
    \begin{subfigure}{0.49\textwidth}
        \centering
        \includegraphics{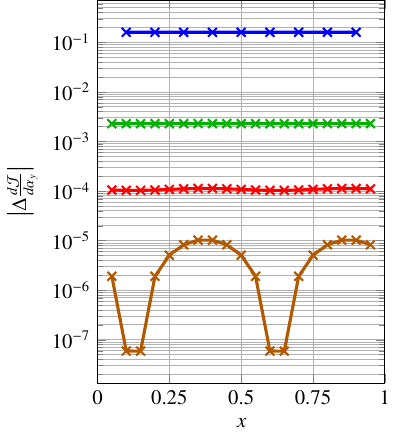}
        \caption{Absolute errors of $\dJday$}
        \label{fig:abs_error_grad_y}
    \end{subfigure}
    \hfill
    \begin{subfigure}{0.49\textwidth}
        \centering
        \includegraphics{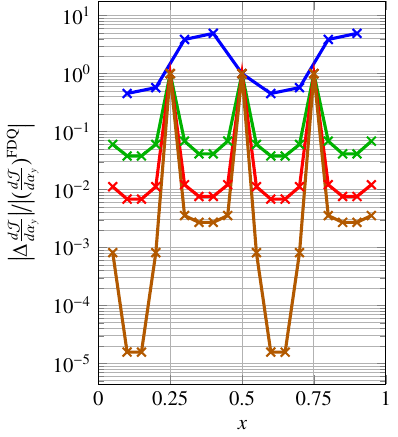}
        \caption{Relative errors of $\dJday$}
        \label{fig:rel_error_grad_y}
    \end{subfigure}
    \vspace{5pt}\\
    \centering \includegraphics{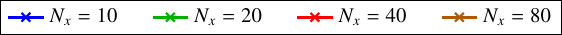}
    \caption{Node-wise validation of the ALBM gradient against the one-sided FDQ
    reference \refEq{eq:FDQ_grad} with $h = 10^{-5}$, evaluated along the diagonal
    $y = x$ for the informed distributed control problem~(ii) of
    \refSec{sec:distributed_periodic} in the first optimization step, with initial
    control  $\al^{(0)}(\x_n) = -\Force^*(\x_n)$.
    Absolute (left column) and relative (right column) deviation of the gradient
    components $\dJdax$ (top row) and $\dJday$ (bottom row), obtained from the
    optimality condition \refEq{eq:optimality}, for resolutions
    $N_x \in \{10, 20, 40, 80\}$. Here
    $\Delta(\bullet) \coloneqq (\bullet)^{\text{ALBM}} - (\bullet)^{\text{FDQ}}$, with both
    gradients evaluated at the same node
    $\x_n = \begin{pmatrix} s & s \end{pmatrix}^{\mathsf{T}}$ of the diagonal. 
    }
    \label{fig:error_grad}
\end{figure}

\subsection{Material parameter identification on a bounded domain \label{sec:param_bounded}}
For the third benchmark, the domain is the elliptic plate of \citet{boolakeeDirichletNeumannBoundary2023} with Dirichlet and Neumann boundaries, and the control is Young's modulus, that is a single scalar. The numerical setup is given in \refSec{sec:numerical_setup_ellipse}, the verification of the primal problem in \refSec{sec:primal_validation_ellipse}, and the identification of the modulus, together with the influence of the initial value, in \refSec{sec:parameter_ellipse}.

\subsubsection{Numerical setup \label{sec:numerical_setup_ellipse}}
\begin{figure}[ht]
    \centering
    \includegraphics{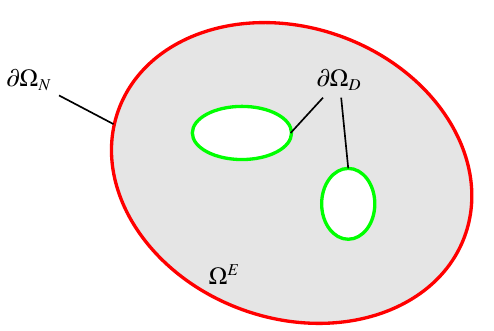}
    \caption{Visualization of the elliptical plate geometry with the respective boundary conditions. Dirichlet conditions (\refEq{eq:disp_Dirichlet}) are indicated in green, and Neumann conditions (\refEq{eq:disp_Neumann}) in red.}
    \label{fig:elliptic_plate_geometry}
\end{figure}

The geometry is the elliptic plate of
\citet{boolakeeDirichletNeumannBoundary2023}, consisting of an outer ellipse from
which two inner ellipses are removed, with the geometric parameters given in
\refTab{tab:ellipseparameter} and the configuration illustrated in
\refFig{fig:elliptic_plate_geometry}. The full numerical domain is
$\Omega = [0, 1.5]^2 \subset \R^2$, discretized by $N_x + 1$ lattice nodes per
Cartesian direction. The elliptic domain on which the PDE constraint
\refEq{eq:Ghat} is enforced is denoted by $\Omega^E$. Its formulation as a level-set function can be found in \cite{boolakeeDirichletNeumannBoundary2023}. The boundary
conditions read
\begin{align}
    \disp &= \disp_D \quad \text{on} \quad \partial \Omega_D \times [0,T],  \label{eq:disp_Dirichlet}\\
    \stress \bm{n} &= \bm{T} \quad \text{on} \quad \partial \Omega_N \times [0,T], \label{eq:disp_Neumann}
\end{align}
where $\disp_D$ is the displacement prescribed on the Dirichlet part
$\partial \Omega_D$ of the boundary, formed by the two inner ellipses, and $\bm{T}$
is the surface traction on the Neumann part $\partial \Omega_N$, formed by the outer
ellipse, with outward unit normal $\bm{n}$, so that
$\partial \Omega^E = \partial \Omega_D \cup \partial \Omega_N$
\cite{boolakeeDirichletNeumannBoundary2023}. The boundary conditions are part of the PDE constraint in \refEq{eq:Ghat}. For the derivation of these boundary conditions on the LBM level and their implementation, the reader is referred to
\cite{boolakeeDirichletNeumannBoundary2023, kaiserFluidStructureInteractionSimulations2025}.

Following \cite{boolakeeDirichletNeumannBoundary2023}, the manufactured steady-state
displacement is given as
\begin{equation}
    \disp^*(\x) = \begin{pmatrix}
        \delta_1 \sin(2\pi x_1 x_2) \\
        \delta_2 \cos(2\pi x_2)(1 + x_1^2)
    \end{pmatrix},
    \label{eq:disp_man_ellipse}
\end{equation}
with $\delta_1 = 9 \cdot 10^{-4}$ and $\delta_2 = 7 \cdot 10^{-4}$, and the
associated stress $\stress^*$ follows from the constitutive relation
\refEq{eq:stress}. The body force is again obtained from the steady Navier--Cauchy
equation,
evaluated for the reference material parameters $E^* = 0.1$ and $\nu = 0.7$. The Dirichlet data $\disp_D$ and the traction vector $\bm{T}$ in
\refEqs{eq:disp_Dirichlet}{eq:disp_Neumann} are taken from
\refEq{eq:disp_man_ellipse} and \refEq{eq:stress} at these reference parameters. The initial
displacement is again $\disp_0(\x, t=0) = \bm{0}$.
\begin{table}[ht]
    \centering
    \begin{tabular}{c c c c c c}
    \hline
        Ellipse & $a_i^2$ & $b_i^2$ & $x_i$ & $y_i$ & $\phi_i$ \\
        \hline
        1 & $0.48L^2$ & $0.3L^2$ & $0.750L$ & $0.750L$ & $-\pi / 9$ \\
        2 & $0.010L^2$ & $0.018L^2$ & $0.950L$ & $0.600L$ & $0$ \\
        3 & $0.035L^2$ & $0.010L^2$ & $0.500L$ & $0.900L$ & $0$ \\
        \hline
    \end{tabular}
    \caption{Geometric parameters of the elliptical plate illustrated in \refFig{fig:elliptic_plate_geometry}, with the characteristic length $L=1$.}
    \label{tab:ellipseparameter}
\end{table}
\subsubsection{Verification of the primal problem \label{sec:primal_validation_ellipse}}
\begin{figure}[!htbp]
    \centering
        \includegraphics{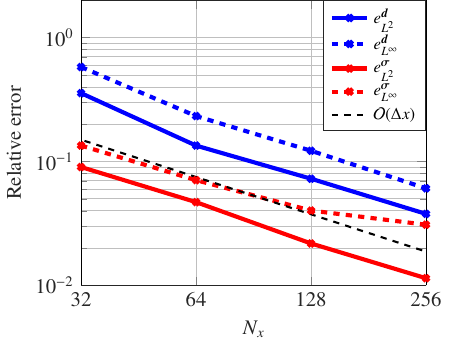}
    \caption{Grid convergence of the relative errors
    $e^{\disp}_{L^p}$ \refEq{eq:err_disp} and $e^{\stress}_{L^p}$ \refEq{eq:err_stress},
    $p \in \{2, \infty\}$, against the manufactured solution
    \refEqs{eq:disp_man_ellipse}{eq:stress}, for
    $N_x \in \{32, 64, 128, 256\}$ lattice nodes per Cartesian direction, evaluated
    after a fixed number of $40 N_x^2$ time steps.}
    \label{fig:primal_BC}
\end{figure}

The primal problem is verified as in \refSec{sec:primal_validation_periodic}, for
$E = 0.1$, $\nu = 0.7$ and $N_x \in \{32, 64, 128, 256\}$. The results are shown in
\refFig{fig:primal_BC}. To exclude the influence of the stopping criterion, every run
is performed with a fixed number of $40 N_x^2$ iterations. Displacement and stress
converge with first order, with $\mathrm{EOC}(e^{\disp}_{L^2}) = 1.06$, $\mathrm{EOC}(e^{\disp}_{L^\infty}) = 1.07$ for the displacement, and $\mathrm{EOC}(e^{\stress}_{L^2}) = 1.01$,  $\mathrm{EOC}(e^{\stress}_{L^\infty}) = 0.71$ 
for the stress, which is one order lower than in the periodic case of
\refFig{fig:prim_periodic}. The reduction is caused by the first-order treatment of the curved Dirichlet and Neumann boundaries, whereas the bulk scheme is of second order
\cite{boolakeeDirichletNeumannBoundary2023, kaiserFluidStructureInteractionSimulations2025}.

\subsubsection{Identification of Young's modulus \label{sec:parameter_ellipse}}
In the third benchmark Young's modulus is the control. To ensure a physically
admissible material parameter, $E > 0$ has to hold, which is enforced here by the
unconstrained substitution
\begin{equation} \label{eq:E_substitution}
    E = \alpha^2 \in \R_{>0}, \qquad \alpha \in \CS_{\Delta x} \subset \R,
\end{equation}
so that the optimizer operates on an unconstrained scalar while the constraint on
$E$ is satisfied by construction. The reference control is
$\alpha^* = \sqrt{E^*}$. The boundary conditions and the exerted body force are held fixed throughout the optimization. These
loads are data of the inverse problem, and only their representation in lattice units
depends on the control. In particular, the
traction $\bm{T} = \stress^*\bm{n}$ has to be evaluated with the reference moduli
$\mu(E^*)$ and $K(E^*)$. If it were evaluated with the control of the respective optimization step, the dependence on $E$ would cancel on $\partial\Omega_N$.

Therefore, the state depends nonlinearly on the control $\alpha$, through the relaxation frequencies of the scheme. At the reference control the
lattice modulus is $\tilde{E}(\alpha^*) = E^* = 0.1$, the value used in the verification
of the primal problem in \refSec{sec:primal_validation_ellipse}.

Gradients are computed with f-AD (\refSec{sec:AD}), which is efficient here because the control is
one-dimensional.
The objective \refEq{eq:discrete_J} is evaluated in the discretized elliptic domain $\Ddesign = \Omega^E_{\Delta x}$. The optimizer parameters are listed in
\refTab{tab:parameter_ellipse}.
\renewcommand{\arraystretch}{1.75}
\begin{table}[H]
    \centering
    \begin{tabular}{|c|c|c|c|c|c|c|c|c|c|}
        \hline
        $T_{max}$ & $\varepsilon_{conv}$ & $\varepsilon$ & $\varepsilon_\alpha$ & $\lambda^0$ & $m_{max}$ & $k_{max}$ & $l$
         & $\gamma^0$\\
        \hline
        $60$ & $10^{-6}$ & $10^{-10}$ &  $2.2 \cdot 10^{-16}$ & 1 & 50 & 100 & 20 & $10^{-4}$\\
        \hline
    \end{tabular}
    \caption{Parameters of the GBO algorithm and the L-BFGS optimizer for the parameter identification problem~(3), with the symbols as defined in
    \refTab{tab:parameter}.}
    \label{tab:parameter_ellipse}
\end{table}
\renewcommand{\arraystretch}{1.0}
A property of the substitution \refEq{eq:E_substitution} is that the map $\alpha \mapsto \alpha^2$ allows for two minimizers $\pm \alpha^*$ and the recovered control is only
determined up to its sign, therefore the error is reported for $|\alpha|$. 

\begin{figure}[!htbp]

    \begin{subfigure}{0.3\textwidth}
    \centering
    \includegraphics{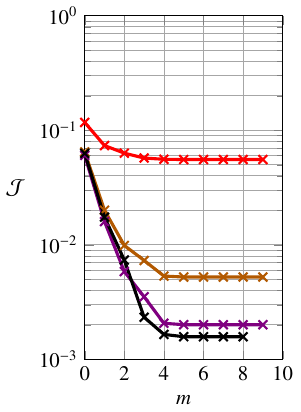}
    \label{fig:obj_param_AD_ellipse}
    \end{subfigure}
    \hfill
    \begin{subfigure}{0.3\textwidth}
    \centering
    \includegraphics{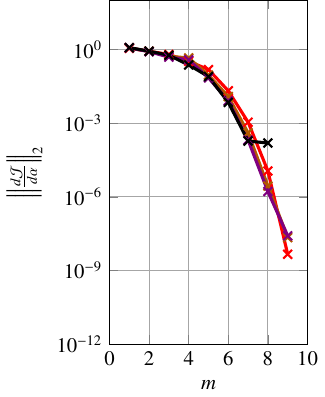}
    \label{fig:deriv_param_AD_ellipse}
    \end{subfigure}
    \hfill
    \begin{subfigure}{0.3\textwidth}
    \centering
    \includegraphics{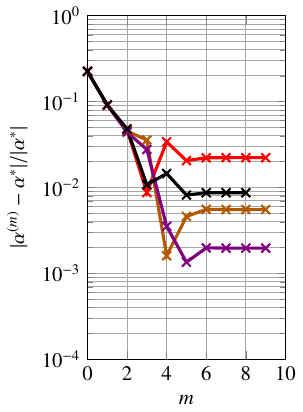}
    \label{fig:controls_param_AD_ellipse}
    \end{subfigure}
    \vspace{5pt}
    \centering \includegraphics{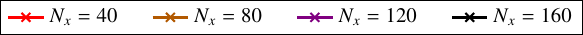}

    \caption{Objective functional $\J$ , magnitude of gradient $\Big|\Big|\dJda \Big| \Big|_2$ and relative error of controls $| \alpha^{(m)} - \alpha^* | / | \alpha^* |$ for parameter identification problem (3) as formulated in \refTab{tab:overview_cases} using f-AD for the gradient computation for resolutions $N_x \in \{ 40,80,120,160 \}$.}
    \label{fig:parameter_AD_ellipse}
\end{figure}

\refFig{fig:parameter_AD_ellipse} shows the evolution of the objective functional
$\J$, of the gradient magnitude and of the relative control error
$|\alpha^{(m)} - \alpha^*|/|\alpha^*|$ over the optimization steps $m$, for
$N_x \in \{40, 80, 120, 160\}$, starting from $E^{(0)} = 0.15$. The optimization
terminates within at most $10$ steps at all resolutions. The termination is triggered
either by the gradient tolerance $\varepsilon$ of \refTab{tab:parameter_ellipse} or,
once the changes of $\J$ between successive line-search trials fall below the level of
the steady-state criterion \refEq{eq:convergence_crit}, by a failure of the line
search. In the latter case the control has settled to at least five digits, so that the
remaining change is far below the discretization error discussed next.

The decisive difference to the periodic benchmarks is the level at which the
objective stagnates. In \refSec{sec:parameter_periodic} and
\refSec{sec:distributed_periodic} the target was exactly representable by the
control and $\J$ decreased to round-off. Here the control is a single scalar that acts on the material parameters,
and no choice of $E$ can compensate the boundary-induced discretization error of the
forward solver. The residual $\disp_{\Delta x}(\alpha) - \disp^*$ therefore remains
finite, and $\J$ stagnates at the level set by the primal error of
\refFig{fig:primal_BC}.

The relative control error stagnates at a resolution-dependent floor as well, with
final values of $2.23\,\%$, $0.55\,\%$, $0.20\,\%$ and $0.87\,\%$ for
$N_x = 40, 80, 120, 160$.
Over the four resolutions the EOC of the control error is $1.06$ \refEq{eq:EOC}, which matches the first-order convergence of the primal problem in \refFig{fig:primal_BC}. Unlike in the periodic benchmarks, however, the error does not decrease monotonically as it falls by an order of magnitude from $N_x = 40$ to $N_x = 120$ and increases again at
$N_x = 160$.
This decrease in the error under refinement is observed for the primal problem in \refFig{fig:primal_BC} between resolutions of $N_x=128$ and $N_x=256$ in the error $e^{\stress}_{L^\infty}$. One possible cause for this behavior could be the way in which curved boundaries cut the lattice.

\begin{figure}[!htbp]

    \begin{subfigure}{0.3\textwidth}
    \centering
    \includegraphics{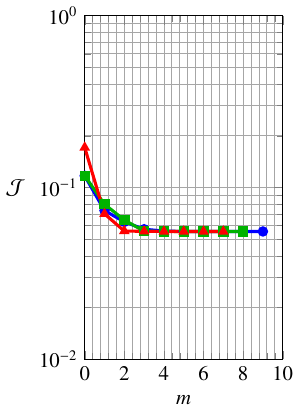}
    \label{fig:obj_param_AD_init_ellipse}
    \end{subfigure}
    \hfill
    \begin{subfigure}{0.3\textwidth}
    \centering
    \includegraphics{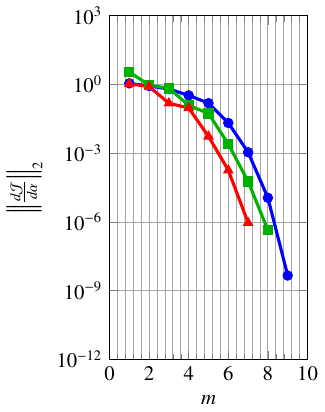}
    \label{fig:deriv_param_AD_init_ellipse}
    \end{subfigure}
    \hfill
    \begin{subfigure}{0.3\textwidth}
    \centering
    \includegraphics{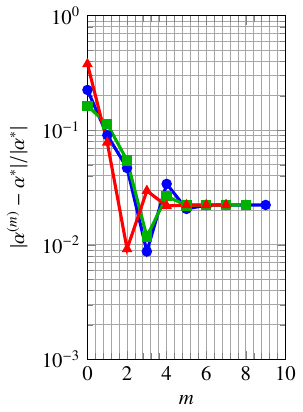}
    \label{fig:controls_param_AD_init_ellipse}
    \end{subfigure}
    \vspace{5pt}
    \centering
    \includegraphics{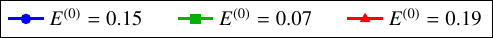}

    \caption{Objective functional $\J$, magnitude of gradient $\Big|\Big|\dJda \Big| \Big|_2$ and relative error of controls $| \alpha^{(m)} - \alpha^* | / | \alpha^* |$ for parameter identification problem (3) as formulated in \refTab{tab:overview_cases} using f-AD for three initial controls $E^{(0)} = [\alpha^{(0)}]^2 \in \{ 0.07, 0.15, 0.19\} $
    for a resolution of $N_x=40$.}
    \label{fig:param_init_controls_ellipse}
\end{figure}

The influence of the initial control is examined in
\refFig{fig:param_init_controls_ellipse}, where three initial values $E^{(0)} = [\alpha^{(0)}]^2 \in \{ 0.07, 0.15, 0.19\} $ are compared at $N_x = 40$. All three runs converge to the same
objective and control-error floor within at most $9$ steps, and the number of
iterations depends only weakly on the initial value ($9$, $8$ and $7$ steps for
$E^{(0)} = 0.15$, $0.07$ and $0.19$). The recovered controls agree to nine digits. This insensitivity is an empirical observation rather than a guaranteed property; it indicates that the objective is well behaved over the range of moduli considered and that the identification is robust with respect to the starting point.

\section{Conclusion} \label{sec:conclusion}
In this work, three benchmark cases for gradient-based optimization in linear
elastic solid mechanics were proposed, in which the PDE constraint is solved
entirely by a lattice Boltzmann scheme. The primal problems are solved by using the MRT schemes of
\citet{boolakeeNewLatticeBoltzmann2023, boolakeeDirichletNeumannBoundary2023} in
their OpenLB implementation \cite{kaiserFluidStructureInteractionSimulations2025},
and the gradients are obtained either by f-AD or
by the discrete ALBM with automatically generated adjoint collision kernels
\cite{itoGenerationEfficientAdjoint2026}. To the best of the authors' knowledge,
this constitutes the first application of the discrete ALBM to PDE-constrained
optimization problems in linear elastic solid mechanics with a purely LBM-based
primal and adjoint solver. Each benchmark provides a manufactured solution
of the forward problem and a known reference control $\al^*$, so that the primal
discretization error, the error of the recovered control and their experimental
orders of convergence are accessible in the same setup.

The first benchmark identifies two scalar amplitudes of an analytical force field
on a periodic domain. Gradients computed with f-AD and with the ALBM agree at every resolution and in every optimization step, so that the
objective, the gradient norm and the control error are indistinguishable between
the two. The two methods are derived by different differentiation methods and share
only the primal solution and the primal collision kernel. Their agreement is therefore a first indicator that the ALBM gradient, and with it the generated adjoint collision kernel, is correct.

The objective decreases to double-precision round-off because the linear scheme with constant coefficients reproduces the single wavenumber of the manufactured force field, so that the solution differs only in the displacement amplitude from the reference solution, which the two scaling parameters compensate exactly. In contrast, the relative control error stagnates at a resolution-dependent floor with an EOC of $1.92$.

The second benchmark expands the control to the entire nodal force field,
$\dim(\al) = N$, which renders f-AD infeasible and leaves the ALBM as the only
practical gradient method. With the control projected onto zero-mean force fields, as
required by the compatibility condition of the periodic problem, both initial controls
converge to round-off in the objective, and the relative error of the recovered force
field attains an EOC of $1.92$ in both cases, the same value as in the two-parameter
benchmark.
The node-wise comparison against a one-sided FDQ stencil confirms the ALBM
gradient along the diagonal of the domain. Away from the zeros of the exact
gradient, both agree to the accuracy of the reference. 

The third benchmark replaces the periodic setting by the bounded elliptic plate
with Dirichlet and Neumann boundaries and identifies Young's modulus as a material parameter. The discretization is fixed independently of the control, so that the modulus enters through the lattice moduli and the relaxation frequencies of the collision, and the state depends nonlinearly on the control. Here the objective does not decrease to round-off but stagnates at the level set by the primal error. 

The relative control error decreases from $2.2\,\%$ at $N_x = 40$ to $0.2\,\%$ at
$N_x = 120$ and rises again to $0.9\,\%$ at $N_x = 160$, with an EOC of $1.06$ over the
four resolutions. It thus follows the first-order convergence of the primal problem on
this domain in magnitude, but not monotonically, which is consistent with the geometry-dependent accuracy of the primal problem on this domain.

The optimization itself proved robust with respect to the initial control in all benchmarks, with all tested initial values reaching the same objective and control-error floor.  
Two limitations exist for the findings of this work. First, adjoint boundary
conditions for the solid LBM scheme were not derived, so the gradient of the
bounded-domain benchmark could only be obtained by f-AD, and the node-wise
validation of the ALBM gradient is restricted to the bulk dynamics. Second, all
cases are two-dimensional and posed in linear elasticity for time-independent problems.

These restrictions define the next steps. The formulation and automatic generation of adjoint boundary conditions would close the remaining gap and make the ALBM
applicable to the bounded-domain benchmark, which is the prerequisite for any
practically relevant structural optimization problem. Beyond that, the framework
extends naturally to non-linear elastic constitutive laws and to shape and topology
optimization of solids. The main motivation, however, is to progress towards time-dependent optimization problems for FSI with a purely LBM-based solver for both phases. The benchmarks proposed here are intended to serve as the verification basis on which such extensions can be built. 

\section*{Funding}
This work has received funding from the European Union’s Horizon Europe research and innovation
program under grant agreement No 101138305.

\section*{Acknowledgements}
The authors thank F. Heberle, E. Herrero and K. Stipcevic who contributed an initial implementation of the bounded-domain optimization setup within a student project. 
The authors thank O. Boolakee who conceptually supported the implementation of the LBM for solids in OpenLB.

\section*{Author contribution statement}
\textbf{J.L.G.:} 
    Conceptualization, 
    Methodology, 
    Software, 
    Validation, 
    Formal Analysis, 
    Investigation, 
    Data Curation, 
    Writing - Original Draft, 
    Writing - Review {\&} Editing, 
    Visualization;
\textbf{F.K.:}
    Methodology, 
    Software, 
    Validation, 
    Data Curation, 
    Writing - Review {\&} Editing, 
    Visualization;
\textbf{S.S.:}
    Methodology,
    Software,
    Validation,
    Formal Analysis,
    Investigation,
    Writing - Review {\&} Editing,
    Supervision,
    Project administration,
    Funding Acquisition;
\textbf{S.I.:} 
    Conceptualization, 
    Methodology, 
    Software, 
    Formal Analysis, 
    Investigation, 
    Writing - Review {\&} Editing, 
    Supervision, 
\textbf{M.J.K.:} 
    Software,
    Writing - Review {\&} Editing, 
    Supervision,
    Resources, 
    Project administration, 
    Funding Acquisition. 
All authors read and approved the final version of this paper.

\section*{Data availability statement}
The results within this work were produced with OpenLB \cite{krauseOpenLBOpenSource2021}, which is released open source under the GNU General Public License, version 2. The computational data are available upon reasonable request. The code used within this work can be found under commit: \href{https://gitlab.com/openlb/olb/-/tree/4a98e2875795fc3bf92ad56319aa0e67a987d257/}{4a98e287} (unreleased).

\section*{Declaration of generative AI in the manuscript preparation process}
During the preparation of this work, the authors used Anthropic Claude and Google Gemini to assist with code development, data analysis, text drafting and formatting. 
After using this tool, the authors reviewed and edited the content as needed and take full responsibility for the content of the publication.

\bibliographystyle{elsarticle-num-names} 
\bibliography{references.bib}

\appendix
\counterwithin*{table}{section}
\counterwithin*{figure}{section}

\section{Solid LBM scheme}
\label{app:solid_LBM}

\subsection{MRT Transformation Matrices}
\label{app:transformation_matrices}
\begin{align}
    \mathbf{M} &=
    \left[\begin{array}{*{8}r}
    1 & 0 & -1 & 0 & \:\:1 & -1 & -1 & 1 \\
    0 & 1 & 0 & -1 & 1 & 1 & -1 & -1 \\
    0 & 0 & 0 & 0 & 1 & -1 & 1 & -1 \\
    1 & 1 & 1 & 1 & 2 & 2 & 2 & 2 \\
    1 & -1 & 1 & -1 & 0 & 0 & 0 & 0 \\
    0 & 0 & 0 & 0 & 1 & -1 & -1 & 1 \\
    0 & 0 & 0 & 0 & 1 & 1 & -1 & -1 \\
    0 & 0 & 0 & 0 & 1 & 1 & 1 & 1
    \end{array}\right],
    \label{app:m_matrix} 
\end{align}
The matrix $\mathbf{M}$ given above is used for the periodic benchmarks~(1) and~(2),
for which $m_{22} = m_f$, i.e.\ $\gamma = 0$ in \refEq{eq:moment_f}. For the
bounded-domain benchmark~(3), the relation $m_{22} = m_f - \gamma m_s$ with
$\tau_f = 1/2$ is used, which replaces the last row of $\mathbf{M}$ by
$(\gamma, \gamma, \gamma, \gamma, 1+2\gamma, 1+2\gamma, 1+2\gamma, 1+2\gamma)$.

\subsection{Relaxation times}
\label{app:relaxation_times}
\renewcommand{\arraystretch}{1.5}
\begin{table}[H]
    \centering
    \begin{tabular}{|c|cccccc|}
         \hline
         & $m_{11}$ & $m_s$ & $m_d$ & $m_{12}$ & $m_{21}$ & $m_f$ \\
         \hline
         $\tau_\beta$ &  $\theta^{-1} \tilde{\mu}$ & $2(1+\theta)^{-1} \tilde{K}$ &
         $2(1-\theta)^{-1} \tilde{\mu}$ &
         $1/2$ & $1/2$ & $1/2$\\
         \hline
    \end{tabular}
    \caption{Relaxation times $\tau_\beta$ for solid LBM scheme, adapted from \cite{boolakeeDirichletNeumannBoundary2023}.}
    \label{tab:relaxation_times}
\end{table}
\renewcommand{\arraystretch}{1}

\end{document}

%% file: commands.tex
\newcommand{\Force}{\bm{F}}
\newcommand{\x}{\bm{x}}
\newcommand{\al}{\bm{\alpha}}
\newcommand{\J}{\mathcal{J}}
\newcommand{\CJ}{\hat{\J}}

\newcommand{\dJda}{\frac{d \mathcal{J}}{d \al}}

\newcommand{\CV}{\mathcal{V}}

\newcommand{\disp}{\bm{d}}
\newcommand{\stress}{\bm{\sigma}}
\newcommand{\normal}{\bm{n}}
\newcommand{\trac}{\bm{T}}

\newcommand{\R}{\mathbb{R}}
\newcommand{\CS}{\mathcal{U}}
\newcommand{\CSparam}{\mathcal{U}_{param}}
\newcommand{\CSdist}{\mathcal{U}_{dist}}
\newcommand{\Cdesign}{\Omega^{\CJ}}

\newcommand{\Ddesign}{\Omega^{\J}_{\Delta x}}

\newcommand{\fracpp}[2]{\frac{\partial #1}{\partial #2}}

\newcommand{\refEq}[1]{(\ref{#1})}
\newcommand{\refEqs}[2]{(\ref{#1}) and (\ref{#2})}
\newcommand{\refFig}[1]{Figure \ref{#1}}
\newcommand{\refTab}[1]{Table \ref{#1}}
\newcommand{\refAlg}[1]{Algorithm \ref{#1}}
\newcommand{\refSec}[1]{Section \ref{#1}}

\newcommand{\fs}{\bm{f}_{\mathcal{S}}}

\newcommand{\fc}{\bm{f}_{\mathcal{C}}}

\newcommand{\fsa}{\bm{\varphi}_{\mathcal{S}}}
\newcommand{\fsai}[1]{\bm{\varphi}_{\mathcal{S}, #1}}
\newcommand{\fca}{\bm{\varphi}_{\mathcal{C}}}
\newcommand{\dCdf}{\fracpp{\mathcal{C}[\fs^\infty, \al]}{\fs}}
\newcommand{\dJdf}{\fracpp{\J(\al, \fs^\infty)}{\fs}}

\newcommand{\dispderivative}{\dJda}
\newcommand{\normdispderivative}{\Big|\Big|\dispderivative \Big|\Big|_{L_2 \left( \Ddesign \right)}}

\newcommand{\relErrorControl}[1]{\left|\left|\bm{F}^{(#1)} - \bm{F^*}\right|\right|_{L^2\left( \Ddesign \right)} / \left|\left|\bm{F^*}\right|\right|_{L^2 \left( \Ddesign \right)}}

\newcommand{\relErrorParameterControl}{\frac{||\bm{\alpha}^{(m)} - \bm{\alpha^*}||_2}{||\bm{\alpha^*}||_2}}